\documentclass[12pt,leqno]{amsart} 
\usepackage[numbers]{natbib}
\usepackage{amsmath,amsthm}
\usepackage{amssymb}
\renewcommand{\widehat}{\hat}
\usepackage{url}
\usepackage[letterpaper,colorlinks={true},
backref,pdfstartview={FitBH},pdfview={FitBH},bookmarks={true}]{hyperref}
\global\hbadness=500
\global\tolerance=200
\newtheorem{theorem}{Theorem}[section]
\newtheorem{lemma}[theorem]{Lemma}
\newtheorem{conjecture}[theorem]{Conjecture}

\newtheorem{question}[theorem]{Question}
\newtheorem{requirement}[theorem]{Requirement}
\newtheorem{corollary}[theorem]{Corollary} 
\theoremstyle{definition}
\newtheorem{definition}[theorem]{Definition}

\theoremstyle{remark}

\newtheorem{remark}[theorem]{Remark}
\numberwithin{equation}{theorem}

\newcommand{\ce}{computably enumerable }
 
\newcommand{\cat}{\widehat{~}}
\newcommand{\join}{\oplus} 
\renewcommand{\phi}{\varphi}

\newcommand{\E}{\mathcal{E}}
\renewcommand{\L}{\mathcal{L}}\renewcommand{\S}{\mathcal{S}}

\newcommand{\F}{\mathcal{F}}
\newcommand{\R}{\mathcal{R}}
\newcommand{\B}{\mathcal{B}}
\newcommand{\D}{\mathcal{D}}

\newcommand{\N}{\mathcal{N}}

\newcommand{\Ahat}{\widehat{A}} 
\newcommand{\Fhat}{\widehat{F}}
\newcommand{\Mhat}{\widehat{M}}
\newcommand{\Dhat}{\widehat{D}}
\newcommand{\Dtilde}{\tilde{D}}
\newcommand{\Rhat}{\widehat{R}}\newcommand{\Rtilde}{\tilde{R}}
\newcommand{\Hhat}{\widehat{H}}
\newcommand{\What}{\widehat{W}}

\newcommand{\That}{\widehat{T}}
\newcommand{\Shat}{\widehat{S}}

\newcommand \cero{\mathbf 0}

\newcommand \Tree{Tree}
\DeclareMathOperator \Isom{Isom}
\newcommand {\w}{\omega}
\newcommand {\andd}{\,\,\&\,\,}
\newcommand {\seq}[1]{\langle #1 \rangle}
\DeclareMathOperator \dom{dom}

\newcommand {\Q}{\mathbb Q}  \newcommand
{\wock}{\omega_1^{\textup{CK}}} 
\newcommand {\HYP}{\textup{HYP}} 

\title[Orbits]{On the Orbits of Computably\\ Enumerable Sets}
  
\author[P.\ Cholak]{Peter~A.~Cholak}

\address{Department of Mathematics\\ University of Notre Dame\\ 
  Notre Dame, IN 46556-5683}

\email{Peter.Cholak.1@nd.edu}

\urladdr{http://www.nd.edu/$\tilde{~}$cholak}

\author[R.\ Downey]{Rodney Downey}

\address{School of Mathematics, Statistics and Computer Science\\
 Victoria University\\ P.O. 
Box 600\\ Wellington, New Zealand}

\email{Rod.Downey@vuw.ac.nz}
\urladdr{http://www.mcs.vuw.ac.nz/$\tilde{~}$downey} 

\author[L.\ Harrington]{Leo~A.~Harrington}

\address{Department of Mathematics\\ University of California \\
  Berkeley, CA 94720-3840}

\email{leo@math.berkeley.edu}

\thanks{Research partially supported NSF Grants DMS-96-34565,
  99-88716, 02-45167 (Cholak), Marsden Fund of New Zealand (Downey),
  DMS-96-22290 and DMS-99-71137 (Harrington).  Some of involved work
  was done partially while Cholak and Downey were visiting the
  Institute for Mathematical Sciences, National University of
  Singapore in 2005. These visits were supported by the Institute.}

\begin{document}

\begin{abstract}
  The goal of this paper is to show there is a single orbit of the
  c.e.\ sets with inclusion, $\E$, such that the question of
  membership in this orbit is $\Sigma^1_1$-complete.  This result and
  proof have a number of nice corollaries: the Scott rank of $\E$ is
  $\wock +1$; not all orbits are elementarily definable; there is no
  arithmetic description of all orbits of $\E$; for all finite $\alpha
  \geq 9$, there is a properly $\Delta^0_\alpha$ orbit (from the
  proof).
\end{abstract}

\maketitle

\section{Introduction}

In this paper we work completely within the c.e.\ sets with inclusion.
This structure is called $\E$.

\begin{definition}
  $A \approx \Ahat$ iff there is a map, $\Phi$, from the c.e.\ sets to
  the c.e.\ sets preserving inclusion, $\subseteq$, (so $\Phi \in
  \text{Aut}(\E)$) such that $\Phi(A)=\Ahat$.
\end{definition}

By \citet{Soare:74}, $\mathcal{E}$ can be replaced with
$\mathcal{E}^*$, $\mathcal{E}$ modulo the filter of finite sets, as
long as $A$ is not finite or cofinite.  The following conjecture was
made by Ted Slaman and Hugh Woodin in 1989.

\begin{conjecture}[\citet{Slaman.Woodin:conjecture}]
  The set $\{\langle i , j \rangle: W_i\approx W_j)\}$ is
  $\Sigma^1_1$-complete.
\end{conjecture}

This conjecture was claimed to be true by the authors in the mid
1990s; but no proof appeared.  One of the roles of this paper is to
correct that omission. The proof we will present is far simpler than
all previous (and hence unpublishable) proofs.  The other important
role is to prove a stronger result.

\begin{theorem}[The Main Theorem]\label{main}
  There is a c.e.\ set $A$ such that the index set $\{i : W_i \approx
  A\}$ is $\Sigma^1_1$-complete.
\end{theorem}
   
As mentioned in the abstract this theorem does have a number of nice
corollaries.

\begin{corollary}
  Not all orbits are elementarily definable; there is no arithmetic
  description of all orbits of $\E$.
\end{corollary}

\begin{corollary}
  The Scott rank of $\E$ is $\wock +1$.
\end{corollary}

\begin{proof}
  Our definition that a structure has Scott rank $\wock +1$ is that
  there is an orbit such that membership in that orbit is
  $\Sigma^1_1$-complete. There are other equivalent definitions of a
  structure having Scott Rank $\wock +1$ and we refer the readers to
  \citet{MR1767842}.
\end{proof}

\begin{theorem}\label{sec:maincor}
  For all finite $\alpha > 8$ there is a properly $\Delta^0_{\alpha}$
  orbit.
\end{theorem}

\begin{proof}
  Section~\ref{sec:last} will focus on this proof.
\end{proof}

\subsection{Why Make Such a Conjecture?}
  
Before we turn to the proof of Theorem~\ref{main}, we will discuss the
background to the Slaman-Woodin Conjecture.  Certainly the set
$\{\langle i , j \rangle: W_i\approx W_j)\}$ is $\Sigma^1_1$.  Why
would we believe it to be $\Sigma_1^1$-complete?

\begin{theorem}[Folklore\footnote{See Section~\ref{folkloreproofI} for
    more information and a proof.}] \label{folkloreBA} There is a
  computable listing, $\mathcal{B}_i$, of computable Boolean algebras
  such that the set $\{\langle i , j \rangle: \mathcal{B}_i\cong
  \mathcal{B}_j\}$ is $\Sigma^1_1$-complete.
\end{theorem}

\begin{definition}
  We define $\L(A)= (\{ W \cup A : W \text{ a c.e.\ set}\}, \subseteq)
  $ and $\L^*(A)$ to be the structure $\L(A)$ modulo the ideal of
  finite sets, $\F$.
\end{definition}

That is, $\L(A)$ is the substructure of $\E$ consisting of all c.e.\
sets containing $A$.  $\L(A)$ is definable in $\E$ with a parameter
for $A$.  A set $X$ is finite iff all subsets of $X$ are
computable. So being finite is also definable in $\E$.  Hence
$\L^*(A)$ is a definable structure in $\E$ with a parameter for $A$.
The following result says that the full complexity of the isomorphism
problem for Boolean algebras of Theorem \ref{folkloreBA} is present in
the supersets of a c.e.\ set.

\begin{theorem}[\citet{Lachlan:68*3}]
  Effectively in $i$ there is a c.e.\ set $H_i$ such that $\L^*(H_i)
  \cong \mathcal{B}_i$.
\end{theorem}

\begin{corollary}
  The set $\{\langle i , j \rangle: \L^*(H_i)\cong \L^*(H_j)\}$ is
  $\Sigma^1_1$-complete.
\end{corollary}
  
Slaman and Woodin's idea was to replace ``$\L^*(H_i)\cong \L^*(H_j)$''
with ``$H_i\approx H_j$''.  This is a great idea which we now know
cannot work, as we discuss below.

\begin{definition}[The sets disjoint from $A$]
  $$\mathcal{D}(A) = ( \{ B : \exists W ( B \subseteq A \cup W \text{ and
  } W \cap A =^* \emptyset)\}, \subseteq).$$ Let
  $\mathcal{E}_{\mathcal{D}(A)}$ be $\mathcal{E}$ modulo
  $\mathcal{D}(A)$.
\end{definition}

\begin{lemma}
  If $A$ is simple then $\mathcal{E}_{\mathcal{D}(A)}
  \cong_{\Delta^0_3} \L^*(A)$.
\end{lemma}

$A$ is \emph{$\mathcal{D}$-hhsimple} iff
$\mathcal{E}_{\mathcal{D}(A)}$ is a Boolean algebra. Except for the
creative sets, until recently all known orbits were orbits of
$\mathcal{D}$-hhsimple sets.  We direct the reader to
\citet{Cholak.Harrington:nd} for a further discussion of this claim
and for an orbit of $\E$ which does not contain any
$\mathcal{D}$-hhsimple sets.  The following are relevant theorems from
\citet{Cholak.Harrington:nd}.

\begin{theorem}\label{dhhsimple}
  If $A$ is $\mathcal{D}$-hhsimple and $A$ and $\Ahat$ are in the same
  orbit then $\E_{\mathcal{D}(A)} \cong_{\Delta^0_3}
  \E_{\mathcal{D}(\Ahat)}$.
\end{theorem}

\begin{theorem}[using \citet{Maass:84}]
  If $A$ is $\mathcal{D}$-hhsimple and simple (i.e., hhsimple) then
  $A\! \approx \Ahat$ iff $\L^*(A) \cong_{\Delta^0_3} \L^*(\Ahat)$.
\end{theorem}

Hence the Slaman-Woodin plan of attack fails.  In fact even more is
true.

\begin{theorem}
  If $A$ and $\Ahat$ are automorphic then
  $\mathcal{E}_{\mathcal{D}(A)}$ and
  $\mathcal{E}_{\mathcal{D}(\Ahat)}$ are $\Delta^0_6$-isomorphic.
\end{theorem}

Hence in order to prove Theorem~\ref{main} we must code everything
into $\D(A)$.  This is completely contrary to all approaches used to
try to prove the Slaman-Woodin Conjecture over the years.  We will
point out two more theorems from \citet{Cholak.Harrington:nd} to show
how far the sets we use for the proof must be from simple sets, in
order to prove Theorem~\ref{main}.

\begin{theorem}
  If $A$ is simple then $A \approx \Ahat$ iff $A \approx_{\Delta^0_6}
  \Ahat$.
\end{theorem}

\begin{theorem}
  If $A$ and $\Ahat$ are both promptly simple then $A \approx \Ahat$
  iff $A \approx_{\Delta^0_3} \Ahat$.
\end{theorem}

\subsection{Past Work and Other Connections}

This current paper is a fourth paper in a series of loosely connected
papers, \citet{mr2003h:03063}, \citet{mr2004f:03077}, and
\citet{Cholak.Harrington:nd}.  We have seen above that results from
\citet{Cholak.Harrington:nd} determine the direction one must take to
prove Theorem~\ref{main}.  The above results from
\citet{Cholak.Harrington:nd} depend heavily on the main result in
\citet{mr2004f:03077} whose proof depends on special
$\mathcal{L}$-patterns and several theorems about them which can be
found in \citet{mr2003h:03063}. It is not necessary to understand any
of the above-mentioned theorems from any of these papers to understand
the proof of Theorem~\ref{main}.

But the proof of Theorem~\ref{main} does depend on Theorems~2.16,
2.17, and 5.10 of \citet{Cholak.Harrington:nd}; see
Section~\ref{sec:oldwork}. The proof of Theorem~\ref{sec:maincor} also
needs Theorem~6.3 of \citet{Cholak.Harrington:nd}. The first two
theorems are straightforward but the third and fourth require
work. The third is what we call an ``extension theorem.''  The fourth
is what we might call a ``restriction theorem''; it restricts the
possibilities for automorphisms.  Fortunately, we are able to use
these four theorems from \citet{Cholak.Harrington:nd} as {\em black
  boxes}.  These four theorems provide a clean interface between the
two papers.  If one wants to understand the {\em proofs} of these four
theorems one must go to \citet{Cholak.Harrington:nd}; otherwise, this
paper is completely independent from its three predecessors.

\subsection{Future Work and Degrees of the Constructed Orbits}

While this work does answer many open questions about the orbits of
c.e.\ sets, there are many questions left open. But perhaps these open
questions are of a more degree-theoretic flavor.  We will list three
questions here.

\begin{question}[Completeness]
  Which c.e.\ sets are automorphic to complete sets?
\end{question}

Of course, by \citet{Harrington.Soare:91}, we know that not every
c.e.\ set is automorphic to a complete set, and partial
classifications of precisely which sets can be found in
\citet{Downey.Stob:92} and Harrington and Soare
\cite{Harrington.Soare:96,Harrington.Soare:98}.

\begin{question}[Cone Avoidance]
  Given an incomplete c.e.\ degree $\mathbf{d}$ and an incomplete
  c.e.\ set $A$, is there an $\Ahat$ automorphic to $A$ such that
  $\mathbf{d} \not\leq_T \Ahat$?
\end{question}

In a technical sense, these may not have a ``reasonable'' answer. Thus
the following seems a reasonable question.

\begin{question}
  Are these arithmetical questions?
\end{question}

In this paper we do not have the space to discuss the import of these
questions.  Furthermore, it not clear how this current work impacts
possible approaches to these questions.  At this point we will just
direct the reader to slides of a presentation of Cholak \cite{2006:c};
perhaps a paper reflecting on these issues will appear later.

One of the issues that will impact all of these questions are which
degrees can be realized in the orbits that we construct in
Theorem~\ref{main} and \ref{sec:maincor}.  A set is \emph{hemimaximal}
iff it is the nontrivial split of a maximal set.  A degree is
\emph{hemimaximal} iff it contains a hemimaximal set.
\citet{Downey.Stob:92} proved that the hemimaximal sets form an
orbit.

We will show that we can construct these orbits to contain at least a
fixed hemimaximal degree (possibly along others) or contain all
hemimaximal degrees (again possibly along others).  However, what is
open is if every such orbit must contain a representative of every
hemimaximal degree or only hemimaximal degrees. For the proofs of
these claims, we direct the reader to Section~\ref{hemimaximal}.

\subsection{Toward the Proof of Theorem~\ref{main}}
The proof of Theorem \ref{main} is quite complex and involves several
ingredients. The proof will be easiest to understand if we introduce
each of the relevant ingredients in context.

The following theorem will prove be to useful.

\begin{theorem}[Folklore\footnote{See Section~\ref{folkloreproofI} for
    more information and a proof.}]\label{folkloreI}
  There is a computable listing $T_i$ of computable infinite branching
  trees and a computable infinite branching tree $T_{\Sigma^1_1}$ such
  that the set $\{ i : T_{\Sigma^1_1} \cong T_i\}$ is
  $\Sigma^1_1$-complete.
\end{theorem}

The idea for the proof of Theorem~\ref{main} is to code each of the
above $T_i$s into the orbit of $A_{T_i}$. Informally let
$\mathcal{T}(A_T)$ denote this encoding; $\mathcal{T}(A_T)$ is defined
in Definition~\ref{sec:TA}.  The game plan is as follows:
  
\begin{enumerate}
\item \textbf{Coding:} For each $T$ build an $A_T$ such that $T \cong
  \mathcal{T}(A_T)$ via an isomorphism $\Lambda\leq_T
  \bf{0}^{(2)}$. (See Remark~\ref{sec:LambdaI} for more details.)
\item \textbf{Coding is preserved under automorphic images:} If $\Ahat
  \approx A_T$ via an automorphism $\Phi$ then $\mathcal{T}(\Ahat)$
  exists and $\mathcal{T}(\Ahat) \cong T$ via an isomorphism
  $\Lambda_{\Phi}$, where $\Lambda_\Phi \leq_T \Phi \join
  \bf{0}^{(2)}$.  (See Lemma~\ref{gameplan2}.)
\item \textbf{Sets coding isomorphic trees belong to the same orbit:}
  If $T \cong \widehat{T}$ via isomorphism $\Lambda$ then $A_T \approx
  A_{\widehat{T}}$ via an automorphism $\Phi_\Lambda$ where
  $\Phi_\Lambda\leq_T \Lambda \join \bf{0}^{(2)}$.
\end{enumerate}

So $A_{T_{\Sigma^1_1}}$ and $A_{T_i}$ are in the same orbit iff
$T_{\Sigma^1_1}$ and $T_i$ are isomorphic.  Since the latter question
is $\Sigma^1_1$-complete so is the former question.

We should also point out that work from \citet{Cholak.Harrington:nd}
plays a large role in part~3 of our game plan; see
Section~\ref{sec:oldwork}.

\subsection{Notation}

Most of our notation is standard.  However, we have two trees involved
in this proof.  We will let $T$ be a computable infinite branching
tree as described above in Theorem~\ref{folkloreI}. For the time being
it will be convenient to think of the construction as occurring for
each tree {\em independently}, but this will later change in Section
\ref{manytrees}.  Trees $T$ we will think of as growing upward. There
will also be the tree of strategies which will denote $Tr$ (which will
grow downward). $\lambda$ is always the empty node (in all trees).
It is standard to use $\alpha,\beta, \delta,\gamma$ to range over
nodes of $Tr$.  We will add the restriction that $\alpha,\beta,
\delta,\gamma$ range \emph{only} over $Tr$.  We will use
$\xi,\zeta,\chi$ to range exclusively over $T$.

\section{The Proof of Theorem~\ref{main}}\label{mainproof}

\subsection{Coding, The First Approximation}

\label{sec:coding1}

The main difficulty in this proof is to build a list of pairwise
disjoint computable sets with certain properties to be described
later.  Work from \citet{Cholak.Harrington:nd}, see
Theorem~\ref{interface}, shows that an essential ingredient to
construct an automorphism between two \ce sets is an extendible
algebra for each of the sets.  In addition, to helping with the
coding, this list of pairwise disjoint computable sets will also
provide the extendible algebras for each of the sets $A_{T_i}$, see
Lemma~\ref{niceextendiblealgebra}.

We are going to assume that we have this list of computable
sets and slowly understand how these undescribed properties arise.
For each node $\chi \in \omega^{<\omega}$ and each $i$, we will build
disjoint computable sets $R_{\chi,i}$.  Inside each $R_{\chi,i}$ we
will construct a c.e.\ set $M_{\chi,i}$.

We need to have an effective listing of these sets.  Fix a computable
one-to-one onto listing $l(e)$ from positive integers to the set of
pairs $(\chi,k)$, where $\chi \in \omega^{<\omega}$ and $k \in
\omega$ such that for all $\chi$ and $n$, if $\xi \preceq \chi$, $m
\leq n$, and $l(i) = (\chi,n)$, then there is a $j \leq i$ such that
$l(j) = (\xi,m)$.  Assume that $l(e) = (\chi,k)$; then we will let
$R_{2e} =R_{\chi,2k}$, $R_{2e+1} =R_{\chi,2k+1}$, $M_{2e} =
M_{\chi,2k}$, and $M_{2e+1} = M_{\chi,2k+1}$.  Which listing of the
$R$s we use will depend on the situation.  We do this as there will be
situations where one listing is evidently better than the other.

\begin{definition}
  $M$ is \emph{maximal} in $R$ iff $M \subset R$, $R$ is a computable
  set, and $M \sqcup \overline{R}$ is maximal.
\end{definition}

The construction will ensure that either $M_{\chi,i}$ will be maximal
in $R_{\chi,i}$ or $M_{\chi,i} =^* R_{\chi,i}$.  If $i$ is odd we will
let $M_{\chi,i} = R_{\chi,i}$.  In this case we say $M_{\chi,i}$ is
\emph{known to be computable}.  This is an artifact of the
construction; the odd sets are errors resulting from the tree
construction.  More details will be provided later.

To build $M_{\chi,i}$ maximal we will use the construction in
Theorem~3.3 of \citet{Soare:87}.  The maximal set construction uses
markers.  The marker $\Gamma_e$ is used to denote the $e$th element of
the complement of the maximal set.  At stage $s$, the marker
$\Gamma_e$ is placed on the $e$th element of the complement of the maximal
set at stage $s$.  In the standard way, we allow the marker $\Gamma_e$
to ``pull'' elements of $\overline{M}_s$ at stage $s+1$ such that the
element marked by $\Gamma_e$ has the highest possible $e$-state and
dump the remaining elements into $M$.

However, at times we will have to destroy this construction of
$M_{\chi,i}$ with some priority $p$.  If we decide that we must
destroy $M_{\chi,i}$ with some priority $p$ at stage $s$ we will just
enumerate the element $\Gamma_p$  is marking into $M_{\chi,i}$ at stage
$s$.  If this occurs infinitely many times then $M_{\chi,i} =^*
R_{\chi,i}$. With this twist, we will just appeal to the construction
in \citet{Soare:87}.

To code $T$, for all $\chi$, such that $\chi \in T$, we will build
pairwise disjoint \ce sets $D_{\chi}$. We will let $A = D_\lambda$.
If $l(i)=(\chi,0)$ then we will let $D_i = D_\chi$. If $l(i) \neq
(\zeta,0)$ then we will let $D_i = \emptyset$. These sets will be
constructed as follows.

\begin{remark}[Splitting $M$]\label{splittingM}
  Let $l(j) = (\chi,i)$. We will use the Friedberg Splitting Theorem;
  we will split $M_{\chi,2i}$ into $i+3$ parts.  Again we will just
  appeal to the standard proof of the Friedberg Splitting Theorem.  We
  will put one of the parts into $D_{\chi}$.  For $0 \leq l \leq i$,
  if $\chi \cat l \in T$ and there is a $j'< j$ such that $l(j') =
  (\chi \cat l, 0)$, then we put one of the parts into $D_{\chi \cat
    l}$.  The remaining part(s) remain(s) disjoint from the union of
  the $D$s; we will name this remaining infinite part $H_{\chi,i}$.
  This construction works even if we later decide to destroy
  $M_{\chi,i}$ by making $M_{\chi,i} =^* R_{\chi,i}$.

  If $M_{\chi,i}$ is known to be computable, we will split
  $R_{\chi,i}$ into $i+3$ \emph{computable} parts distributed as
  above.  However in this case we cannot appeal to the Friedberg
  Splitting Theorem since many of the elements in the $D$ under
  question will have entered the $D$s prior to entering $M_{\chi,i} =
  R_{\chi,i}$.  We will have to deal with this case in more detail
  later.
\end{remark}

\begin{lemma}
  This construction implies that $\bigsqcup_\chi D_\chi \subseteq
  \bigsqcup_{(\chi,i)} (R_{\chi,i} - H_{\chi,i})$.
\end{lemma}

At this point we should point out  a possible problem.  If the list of
computable sets is effective then we have legally constructed c.e.\
sets. If not, we could be in trouble.

However, we want our list to satisfy the following requirement.  This
requirement will have a number of roles.  Its main function is to
control where the sets $W_e$ live within our construction.

\begin{requirement}\label{sec:requireI}
  For all $e$, there is an $i_e$ such that either
  \begin{equation}
    \label{eq:1}
    W_e  \cup \bigsqcup_{j \leq i_e} R_i \cup
    \bigsqcup_{j \leq i_e} D_i =^* \omega \text{, or}
  \end{equation}
  \begin{equation}
    \label{eq:2}
    W_e \subseteq^*
    \bigsqcup_{j \leq i_e} R_i  \text{, or}
  \end{equation}
  \begin{equation}
    \label{eq:3}
    W_e \subseteq^*
    \bigsqcup_{j \leq i_e} R_i \sqcup \biggl ( \bigsqcup_{j<i_e} D_i 
    - \bigsqcup_{j \leq i_e} R_i \biggr ).
  \end{equation}
\end{requirement}

Equation~(\ref{eq:3}) implies Equation~(\ref{eq:2}), but this
separation will be useful later.  If Equation~\ref{eq:1} holds, then
there is a computable $R_{W_e}$ such that
\begin{equation}
  \label{eq:8}
  R_{W_e} \subseteq
  \bigsqcup_{j \leq i_e} R_i \cup \bigsqcup_{j \leq i_e} D_i \text{ and }
  W_e \cup R_{W_e} = \omega.
\end{equation}

If we have an effective list of all the $R_{e}$ then we have an
effective list of $H_{e}$. Let $h_i$ be the $i$th element of $H_i$.
Then the collection of all $h_i$ is a computable set, say $W_e$. But
$e$ contradicts Requirement~\ref{sec:requireI}.  It follows that our
list {\em cannot} be effective, but it will be effective enough to
ensure the $D$ are computably enumerable.

At this point we are going to have to bite the bullet and admit that
there will be an underlying tree construction.  We are going to have
to decide how the sets we want to construct will be placed on the
tree.

Assume that $\alpha$ is in our tree of strategies and $l(|\alpha|) =
(\chi,n)$.  At node $\alpha$ we will construct two computable sets
$R_\alpha$ and $E_\alpha$. $E_\alpha$ will be the error forced on us
by the tree construction.  If $\chi \in T$ and $n = 0$ then at
$\alpha$ we will also construct $D_\alpha$.

Assume $\alpha$ is on the true path and $l(|\alpha|) = (\chi,n)$. Then
$R_{\chi,2n} = R_\alpha$ and $E_{\alpha}$ is $R_{\chi,2n+1} =
M_{\chi,2n+1} = E_\alpha$. This is the explanation of why $M_{\chi,i}$
is computable for $i$ odd; $R_{\chi,i}$ is the error.  If $\chi \in T$
and $n = 0$ then $D_\chi = D_\alpha$.  Hence the listing of computable
sets we want is along the true path.  Therefore, from now on, when we
mention $R_{\chi,i}$, $D_\chi$, $R_{e}$, or $D_e$, we assume we are
working along the true path.  When we mention $R_\alpha$ or $D_\alpha$
we are working somewhere within the tree of strategies but not
necessarily on the true path.

\subsection{Meeting Requirement~\ref{sec:requireI}}
 
Our tree of strategies will be a $\Delta^0_3$ branching tree. Hence at
$\alpha$ we can receive a guess to a finite number of $\Delta^0_3$
questions asked at $\alpha^-$.  Using the Recursion Theorem these
questions might involve the sets $R_\beta, E_\beta$, and $D_\beta$ for
$\beta \prec \alpha$.  The correct answers are given along the true
path, $f$.  There is a standard approximation to the true path, $f_s$.
Constructions of this sort are found all over the c.e.\ set
literature.

These constructions are equipped with a computable position function
$\alpha(x,s)$, the node in $Tr$ where $x$ is at stage $s$. All balls
$x$ enter $Tr$ at $\lambda$.  If the approximation to the true path is
the left of $x$'s position, $x$ will be moved upward to be on this
approximation and never allowed to move right of this approximation.
To move a ball $x$ downward from $\alpha^-$ to $\alpha$, $\alpha$ must
be on the approximation to the true path and $x$ must be $\alpha^-$
allowed.  When we $\alpha^-$ allow $x$ depends on Equations~\ref{eq:1}
and \ref{eq:3}.

So, formally, $\alpha(x,x) = \lambda$. If $f_{s+1} <_L \alpha(x,s)$
then we will let $\alpha(x,s+1) = f_{s+1} \cap \alpha(x,s)$.  If
$\alpha(x,s) = \alpha^-$, $x$ has been \emph{$\alpha^-$ allowed},
$\alpha \subseteq f_s$, and, for all stages $t$, if $x \leq t < s$
then $f_t \not<_L \alpha$; then we will let $\alpha(x,s+1) = \alpha$.

Exactly when a ball will be $\alpha$-allowed is the key to this
construction and will be addressed shortly. However, given these rules,
it is clear if $f <_L \alpha$ then there are no balls $x$ with $\lim_s
\alpha(x,s) = \alpha$ and if $\alpha <_L f_s$ then there are at most
finitely many balls $x$ with $\lim_s \alpha(x,s) = \alpha$. Of course,
the question remains what happens at $\alpha \subset f$?

The question we ask at $\alpha^-$ is if the set of $x$ such that there
is a stage $s$ with
\begin{equation}
  \label{eq:5}
  \begin{split}
    x \in W_{e,s}, \alpha^- \subseteq \alpha(x,s), x \text{ is }
    \alpha^- \text{-allowed at stage $s$,} \\
    \text{and } x \notin (\bigsqcup_{\beta \preceq \alpha^-}
    R_{\beta,s} \cup \bigsqcup_{\beta \preceq \alpha^-} E_{\beta,s}
    \cup \bigsqcup_{\beta \preceq \alpha^-} D_{\beta,s})
  \end{split}
\end{equation}
is infinite, where $e = |\alpha^-|$, a $\Pi^0_2$ question.

\subsubsection{A Positive Answer}

Assume that $\alpha$ believes the answer is yes. Then for each time
$\alpha \subset f_s$, $\alpha$ will be allowed to pull three such
balls to $\alpha$.  That is, $\alpha$ will look for three balls
$x_1,x_2,x_3$ and stages $t_1, t_2,t_3$ such that Equation~\ref{eq:5}
holds for $x_i$ and $t_i$, $x_i > s$, $\alpha(x_i,t_i) \not<_L
\alpha$, $x_i \not\in E_{\alpha,t_i} \cup R_{\alpha,t_i}$, and $x$ is
not $\alpha$-allowed at stage $t_i$.

When such a ball $x_i$ and stage $t_i$ are found, we will let
$\alpha(x_i,t_i+1) = \alpha$.  For the first such ball $x_1$ we will
add $x_1$ to $E_\alpha$ at stage $t_1$. Throughout the whole stagewise
construction we will enumerate $x_1$ into various disjoint $D_\beta$
at stage $t_1$ to ensure that $H_\alpha = E_\alpha - \bigsqcup_{\beta
  \preceq \alpha} D_\beta$ and, for each $\beta \preceq \alpha$,
$D_\beta \cap E_\alpha$ is an infinite set. For the second such ball
$x_2$ we will add $x_2$ to $R_\alpha$ at stage $t_2$.  For the third
such ball $x_3$ we will $\alpha$-allow $x_3$ and place all balls $y$
such that $\alpha(y,t_3) = \alpha$, $y \not\in R_{\alpha,t_3}$, and
$y$ is not $\alpha$-allowed into $E_{\alpha,t_3}$ (without any extra
enumeration into the $D_\beta$).

It is not hard to see that when balls are $\alpha$-allowed at stage
$s$ they are not in 
\[
\bigsqcup_{\beta \preceq \alpha^-} R_{\beta,s}
\cup \bigsqcup_{\beta \preceq \alpha^-} E_{\beta,s} \cup
\bigsqcup_{\beta \preceq \alpha} D_{\beta,s};
\]
 once a ball is
$\alpha$-allowed it never enters $R_\alpha$ or $E_\alpha$, and, for
almost all $x$, if $\lim_s \alpha(x,s) = \alpha$ then $x \in E_\alpha
\sqcup R_\alpha$ (finitely many of the $\alpha$-allowed balls may live
at $\alpha$ in the limit).

Assume $\alpha \subset f$.  Then every search for a triple of such
balls will be successful; both $R_\alpha$ and $E_\alpha$ are disjoint
infinite computable sets; infinitely many balls are $\alpha$-allowed
and hence almost of the $\alpha$-allowed balls move downward in $Tr$;
$E_\alpha - \bigsqcup_{\beta \preceq \alpha} D_\beta$ is infinite and
computable; for each $\beta \preceq \alpha$, $D_\beta \cap E_\alpha$
is infinite and computable; $R_\alpha \subset W_e$, and most
importantly, for all $\beta \succ \alpha$, $R_\beta \sqcup
E_\beta\subseteq W_e$ and hence Equation~\ref{eq:1} holds.

\subsubsection{A Negative Answer}

Assume that $\alpha$ believes the answer is no. Assume $\alpha \subset
f$ and that infinitely many balls are $\alpha^-$ allowed. This is 
certainly the case if $\alpha^-$ corresponds to the above positive
answer.  If $W_e$ intersect the sets of balls which are
$\alpha^-$-allowed is finite then 
\[
W_e \subseteq^* \bigsqcup_{\beta
  \preceq \alpha^-} R_{\beta} \cup \bigsqcup_{\beta \preceq \alpha^-}
E_{\beta}
\]
 and hence Equation~\ref{eq:2} holds.  Assume this is not
the case. Since Equation~\ref{eq:5} does not hold for infinitely many
balls $x$ and stages $s$, for almost all $x$ if
\[
x \in W_{e,s}, \alpha^- \subseteq \alpha(x,s), x \text{ is }
\alpha^- \text{-allowed at stage $s$,}
\]
 then $x \in \bigsqcup_{\beta
  \preceq \alpha^-} D_{\beta,s}$.  Hence, 
  \[
  W_e \subseteq^*
\bigsqcup_{\beta \preceq \alpha^-} R_{\beta} \cup \bigsqcup_{\beta
  \preceq \alpha^-} E_{\beta} \cup \bigsqcup_{\beta \preceq \alpha^-}
D_{\beta}
\]
 and Equation~\ref{eq:3} holds.

Either way there are infinitely many balls $x$ and stage $s$ such that
\begin{equation}
  \label{eq:6}
  \begin{split}
    \alpha^- \subseteq \alpha(x,s), x \text{ is }
    \alpha^- \text{-allowed at stage $s$,} \\
    \text{and } x \notin (\bigsqcup_{\beta \preceq \alpha^-}
    R_{\beta,s} \cup \bigsqcup_{\beta \preceq \alpha^-} E_{\beta,s}
    \cup \bigsqcup_{\beta \preceq \alpha^-} D_{\beta,s})\,.
  \end{split}
\end{equation}
In the same way as when $\alpha$ corresponds to the positive answer,
we will pull three such balls to $\alpha$.  The action we take with
these balls is exactly the same as in the positive answer.  Hence,
among other things, infinitely many balls are $\alpha$-allowed, allowing us
to inductively continue.

\subsubsection{The maximal sets and their
  splits}\label{sec:maximal-sets-their}

To build $M_\alpha$ we will appeal to the standard maximal set
construction as suggested above.  But we will label the markers as
$\Gamma^\alpha_e$ or $\Gamma^{\chi,i}_e$ rather than $\Gamma_e$ just
to keep track of things.  As suggested in Remark~\ref{splittingM}, to
build the $D_\beta$ within $R_\alpha$, for $\beta \preceq \alpha$, we
will appeal to the Friedberg Splitting Theorem.

At this point, we will step away from the construction and see what we
have manged to achieve and what more needs to be achieved. We will be
careful to point out where we use the above requirement and where it
is not enough for our goals.

\subsection{A definable view of our coding}

For each $\chi \in T$ we will construct pairwise disjoint c.e.\ sets
$D_\chi$.  The reader might wonder how this helps.  In particular, how
do these sets code $T$?  Moreover, if $\Ahat$ is in the orbit of $A$
how do we recover an isomorphic copy of $T$?  To address these issues,
we will need some sort of ``definable structure.''  Unfortunately, the
definition of the kind of structure we need is rather involved.  To
motivate the definition, we need to recall how nontrivial splits of
maximal sets behave and then see what the above construction does with
these splits in a definable fashion.

\begin{definition}
  A split $D$ of $M$ is a \emph{Friedberg split} iff, for all $W$, if
  $W - M$ is not a c.e.\ set then neither is $W-D$.
\end{definition}

\begin{lemma}[\citet{Downey.Stob:92}]\label{sec:friedbergsplit}
  Assume $M$ is maximal in $R$. Then $D$ is a nontrivial split of $M$
  iff $D$ is a Friedberg split of $M$.
\end{lemma}

\begin{proof}
  In each direction we prove the counterpositive.  Let $\breve{D}$ be
  such that $D \sqcup \breve{D} = M$.

  Assume that $D$ is not Friedberg.  Hence for some $W$, $W - D$ is
  c.e.\ but $W - M$ is not.  If $W \subseteq^* (M \cup \overline{R})$
  then $(W-M)\subseteq^* \overline{R}$ and hence $W-M =^* W \cap
  \overline{R}$, a c.e.\ set.  Therefore $ \overline{M \sqcup
    \overline{R}} = (R-M) \subseteq^* W$.  Therefore $D \sqcup (( W -
  D) \cup \breve{D} \cup \overline{R}) = \omega$ and $D$ is
  computable.

  The set $R - M$ is not a c.e.\ set.  Assume $D$ is computable.  Then
  $R - D = R \cap \overline{D}$.  Hence $\overline{D}$ witnesses that
  $D$ is not a Friedberg split.
\end{proof}

\begin{lemma}\label{sec:UM}
  Assume that $M_i$ are maximal in $R$ and $D$ is a nontrivial split
  of both $M_i$.  Then $M_1 =^* M_2$.
\end{lemma}

\begin{proof}
  $M_1 \cup \overline{R}$ is maximal. $\overline{M_1 \cup
    \overline{R}} = R- M_1$.  Since $M_2\cup \overline{R}$ is maximal
  either $M_1 \subseteq^* M_2$ or $(R -M_1) \subseteq^* M_2$.  In the
  former case, $M_2 \subseteq^* M_1 \cup \overline{R}$ so $M_1 =^*
  M_2$.

  Assume the later case. Let $D \cup \breve{D} = M_2$.  Since $D$ is a
  split of $M_1$, $(R -M_1) \subseteq^* \breve{D}$.  Now $\breve{D} -
  M_1 = R - M_1$ is not c.e.\ set but $\breve{D} - D = \breve{D}$ is a
  c.e.\ set.  So $D$ is not a Friedberg split of $M_1$.  So by
  Lemma~\ref{sec:friedbergsplit}, $D$ is not nontrivial split of
  $M_1$. Contradiction.
\end{proof}

It turns out that we will need a more complex version of the above
lemmas.

\begin{definition}
  $W \equiv_{\mathcal{R}} \What $ iff $W \triangle \What = (W - \What)
  \sqcup (\What - W)$ is computable.
\end{definition}

\begin{lemma}\label{sec:livesI}
  Assume that $M_i$ is maximal in $R_i$ and $D \cap R_i$ is a
  nontrivial split of $M_i$.  Either
  \begin{enumerate}
  \item there are disjoint $\Rtilde_i$ such that $(M_i \cap
    \Rtilde_i)$ is maximal in $\Rtilde_i$, $D \cap \Rtilde_i$ is a
    nontrivial split of $M_i$, and either $\Rtilde_1 = R_1 - R_2$ and
    $\Rtilde_2 =R_2$ or $\Rtilde_1 =R_1$ and $\Rtilde_2 = R_2 - R_1$,
    or
  \item $\tilde{M} = M_1 \cap M_2$ is maximal in $\Rtilde = R_1 \cap
    R_2$.  So $\tilde{R} - M_i =^* \tilde{R} - \tilde{M}$ and hence
    $\tilde{M} \equiv_\R M_1 \equiv_\R M_2$.  Furthermore, if $R_1 =
    R_2$ then $\tilde{M} =^* M_1 =^* M_2$.
  \end{enumerate}
\end{lemma}

\begin{proof}
  $M_i \cup \overline{R}_i$ is maximal. $\overline{M_i \cup
    \overline{R_i}} = R_i- M_i$.  $R_i - M_i$ is not split into two
  infinite pieces by any c.e.\ set.  Since $M_2\cup \overline{R}$ is
  maximal either $(M_1 \cup \overline{R}_1) \subseteq^* (M_2 \cup
  \overline{R}_2)$ or $(R_1 - M_1) \subseteq^* (M_2 \cup
  \overline{R}_2)$.  If $(R_1 - M_1) \subseteq^* (M_2 \cup
  \overline{R}_2)$ then $(R_1 - M_1) \subseteq^* M_2$ or $(R_1 - M_1)
  \subseteq^* \overline{R}_2$.

  Assume $(R_1 - M_1) \subseteq^* M_2$.  So $M_2 - (M_1 \cup
  \overline{R}_1) = R_1 - M_1$ is not a c.e.\ set. Let $(D \cap R_2)
  \cup \Dtilde = M_2$.  Therefore $(R_1 -M_1) \subseteq^* \Dtilde$ or
  $(R_1 -M_1) \subseteq^* (D \cap R_2)$.  In the former case $(D\cap
  R_2) - (M_1 \cup \overline{R}_1) = \emptyset$ is a c.e.\ set. In the
  latter case $\Dtilde - (M_1 \cup \overline{R}_1) = \emptyset$ is a
  c.e.\ set.  Either way, by Lemma~\ref{sec:friedbergsplit}, $(D\cap
  R_2)$ is not a nontrivial split of $M_2$. Contradiction.

  Now assume $(R_1 - M_1) \subseteq^* \overline{R}_2$. Let $\Rtilde_1
  = R_1 - R_2$ and $\Rtilde_2 = R_2$.  Let $(D \cap R_1) \sqcup
  \Dtilde = M_1$ be a nontrivial split.  Let $\tilde{M} = M_1 - R_2$.
  Then $(D \cap \Rtilde_1) \sqcup (\tilde{D} - R_2) = \tilde{M}$ is a
  nontrivial split of $\tilde{M}$.  (Otherwise $(D \cap R_1) \sqcup
  \Dtilde = M_1$ is a trivial split.)

  We can argue dually switching the roles of $M_1$ and $M_2$.  We are
  left with the case $(M_1 \cup \overline{R}_1) \subseteq^* (M_2 \cup
  \overline{R}_2)$ and $(M_2 \cup \overline{R}_2) \subseteq^* (M_1
  \cup \overline{R}_1)$.  Hence $(M_1 \cup \overline{R}_1) =^* (M_2
  \cup \overline{R}_2)$ and $R_1 - M_1 =^* R_2 -M_2$.  Therefore
  $\tilde{M} = M_1 \cap M_2$ is maximal in $\Rtilde = R_1 \cap R_2$.
\end{proof}

\begin{definition}\label{sec:livesinside}
  $D$ \emph{lives inside} $R$ \emph{witnessed by $M$} iff $M$ maximal
  in $R$ and $D \cap R$ is a nontrivial split of $M$.
\end{definition}
 
By Lemma~\ref{sec:UM}, if $D$ lives in $R$ witnessed by $M_i$ then
$M_1 =^* M_2$. Hence at times we will drop the ``witnessed by $M$.''
If $D$ lives in $R$ then we will say $D$ lives in $R$ witnessed by
$M^R$.  The point is that $M^R$ is well defined modulo finite
difference.

\begin{lemma}\label{sec:definable-view-our-2}
  If $D$ lives in $R_1$, $R_1 \cap R_2 = \emptyset$, and $D \cap R_2$
  is computable, then $D$ lives in $R_1 \sqcup R_2$.
\end{lemma}

\begin{lemma}\label{sec:definable-view-our-4}
  If $R$ is computable and $D \cap R$ is computable, then $D$ does not
  live in $R$.
\end{lemma}

\begin{lemma}\label{sec:definable-view-our-3}
  If $\chi \in T$, then $D_\chi$ lives in $R_{\chi,2i}$ or $M_{\chi,i}
  =^* R_{\chi,i}$.
\end{lemma}

\begin{proof}
  Follows from the construction.
\end{proof}

\begin{lemma}\label{sec:definable-view-our-8}
  For all $R_{\chi,i}$, if $M_{\chi,i}$ is maximal in $R_{\chi,i}$,
  there is a subset $H_{\chi,i} \subset M_{\chi,i}$ such that
  $H_{\chi,i}$ lives in $R_{\chi,i}$ and $H_{\chi,i} \cap
  \bigsqcup_{\xi} D_{\xi} = \emptyset$.
\end{lemma}

\begin{proof}
  Follows from the construction.
\end{proof}

\begin{lemma}\label{sec:definable-view-our-5}
  If $D_\xi \cap R_{\chi,i} \neq \emptyset$, then $\xi = \chi$ or
  $|\xi| = |\chi|+1$.  Furthermore, if $D_\xi$ lives in $R_{\chi,i}$
  then $i$ is even.
\end{lemma}

\begin{proof}
  Again follows from the construction.
\end{proof}

\begin{lemma} \label{sec:definable-view-our} If $\chi \cat l \in T$
  then there is a least $i'$ and $j'$ such that $l(j') = (\chi,i')$,
  and, for all $i \geq 2i'$, $D_\chi \cap R_{\chi,i} \neq^*
  \emptyset$, $D_{\chi \cat l}\cap R_{\chi,i} \neq^* \emptyset$, and
  either both $D_\chi$ and $D_{\chi \cat l}$ live in $R_{\chi,i}$ or
  $M_{\chi,i} =^* R_{\chi,i}$.  So, in particular, both $D_\chi$ and
  $D_{\chi \cat l}$ live in $R_{2j'}$ or $M_{2j'} =^* R_{2j'}$.
  Furthermore $i'$ and $j'$ can be found effectively.
\end{lemma}

\begin{proof}
  Assume $\chi \cat l \in T$.  Let $j$ be such that $l(j) = (\chi\cat
  l,0)$.  Let $j'$ be the least such that $j < j'$ and $l(j') =
  (\chi,i')$.  (See Section~\ref{sec:maximal-sets-their}.)
\end{proof}

\begin{requirement}\label{require3}
  For each $\chi \in T$ there are infinitely many $i$ such that
  $M_{\chi,i} \neq^* R_{\chi,i}$.
\end{requirement}

Currently we meet this requirement since if $i$ is even then
$M_{\chi,i} \neq^* R_{\chi,i}$.  But for later requirements we will
have to destroy some of these $M_{\chi,i}$, so some care will be
needed to ensure that it is met.

The following definition is a complex inductive one.  This definition
is designed so that if $A$ and $\Ahat$ are in the same orbit witnessed
by $\Phi$ we can recover a possible image for $D_\chi$ without knowing
$\Phi$.  In reality, we want more: we want to be able to recover
$T$. But the ability to recover $T$ will take a lot more work.  In any
case, the definition below is only a piece of what is needed.

\begin{definition}\label{sec:definable-view-our-6}\mbox{}
  \begin{enumerate}
  \item \emph{An $\mathcal{R}^{A}$ list} (or, equivalently, an
    $\mathcal{R}^{D_\lambda}$ list) is an infinite list of disjoint
    computable sets $R^A_i$ such that, for all $i$, $A$ lives in
    $R^A_i$ witnessed by $M^A_i$ and, for all computable $R$, if $A$
    lives in $R$ witnessed by $M$ then there is exactly one $i$ such
    that $R-M =^* R^A_i-M^A_i$.
  \item We say that \emph{$D$ is a $1$-successor of $\Dtilde$ over
      some $\mathcal{R}^{\Dtilde}$ list} if $D$ and $\Dtilde$ are
    disjoint, and, for almost all $i$, $D$ lives in $R^{\Dtilde}_i$.
  \item Let $D$ be a $1$-successor of $\Dtilde$ witnessed by an
    $\mathcal{R}^{\Dtilde}$ list.  \emph{An $\mathcal{R}^{D}$ list
      over an $\mathcal{R}^{\Dtilde}$ list} is an infinite list of
    disjoint computable sets $R^D_i$ such that, for all $i$, $D$ lives
    in $R^D_i$ and, for all computable $R$, if $D$ lives in $R$ then
    there is exactly one $i$ such that exactly one of $R-M =^* R^D_i-
    M^D_i$ or $R-M =^* R^{\Dtilde}_i-M^{\Dtilde}_i$ hold.
  \end{enumerate}
\end{definition}

\begin{lemma}\label{sec:definable-view-our-7}
  If $\chi \in T$, then let $R^{D_\chi}_e = R_{\chi,g(e)}$, where
  $g(e)$ is the $e$th set of all those $R_{\chi,i}$ where $M_{\chi,i}
  \neq^* R_{\chi,i}$. (By Requirement~\ref{require3}, such a $g$
  exists).  This list is an $\mathcal{R}^{D_\chi}$ list over
  $\mathcal{R}^{D_{\chi^-}}$ (where $\R^{D_\lambda^-}$ is the empty
  list.)
\end{lemma}

\begin{proof}
  We argue inductively.  We are going to take two lists
  $\mathcal{R}^{D_{\chi^-}}$ and $\mathcal{R}^{D_\chi}$ and merge them
  to get a new list.  To each set of this new list we will add at most
  finitely different $R_{\xi,j}$, where for all $i$, $R_{\xi,j}
  -M_{\xi,j} \neq^* R^{D_{\chi^-}}_i-M^{D_{\chi^-}}_i$ and $R_{\xi,j}
  -M_{\xi,j} \neq^* R^{D_{\chi}}_i-M^{D_{\chi}}_i$ such that all such
  $R_{\xi,j}$ are added to some set in our new list.  Call the $n$th
  set of this resulting list $\Rtilde_n$. By
  Lemmas~\ref{sec:definable-view-our-3} and
  \ref{sec:definable-view-our-2} and
  Definition~\ref{sec:definable-view-our-6}, $D_\chi$ lives in almost
  all $\Rtilde_n$.

  Fix $R$ such that $D_\chi$ lives in $R$. For each $n$, apply
  Lemma~\ref{sec:livesI} to $R$ and $\Rtilde_n$. If case (2) applies,
  then $R$ behaves like $\Rtilde_n$ and we are done.  Otherwise we can
  assume $R$ is disjoint from $\Rtilde_n$.

  If this happens for all $n$ then $R$ and $\bigsqcup_i \Rtilde_i$ are
  disjoint.  Split $R$ into two infinite computable pieces $R_1$ and
  $R_2$. Since $\bigsqcup D \subseteq \bigsqcup \Rtilde$, $R_i$ cannot
  be a subset of $\bigsqcup D$.  Therefore $R_i \not\subseteq^*
  \bigsqcup \Rtilde \cup \bigsqcup D$.  Furthermore, $R_i \cup
  \bigsqcup \Rtilde \cup \bigsqcup D \neq^* \omega$.  But assuming
  that we meet Requirement~\ref{sec:requireI} this cannot occur.
  Contradiction.
\end{proof}

\begin{corollary}
  Assume $\chi \cat l \in T$.  By
  Lemmas~\ref{sec:definable-view-our-7} and
  \ref{sec:definable-view-our}, $D_{\chi \cat l}$ is a $1$-successor
  of $D_\chi$ over $\R^{D_\chi}$.  Furthermore, if $F$ is finite then
  $D_{\chi \cat l} - \bigsqcup_{i \in F}R_i$ is a $1$-successor of
  $D_\chi$ over $\R^{D_\chi}$.
\end{corollary}

\begin{corollary}
  If disjoint $D_i$ are $1$-successors of $\Dtilde$ over
  $\R^{\Dtilde}$ then so is $D_1 \sqcup D_2$. In particular, for all
  $\chi, \zeta \in T$, if $\chi \neq \zeta$ and $|\chi| = |\xi|$, then
  $D= D_\chi \sqcup D_\xi$ is a $1$-successor of $D_{\chi^-}$ over
  $\R^{D_{\chi^-}}$ and the elementwise union of the lists
  $\R^{D_\chi}$ and $\R^{D_\xi}$ is an $\R^{D}$ list over
  $\R^{D_{\chi^-}}$.
\end{corollary}

\begin{lemma}
  If $\chi$ does not have a successor in $T$ then there are no
  $1$-successors of $D_\chi$ over $\R^{D_\chi}$.
\end{lemma}

\begin{proof}
  Assume that $D$ is a $1$-successor of $D_\chi$ over
  $\R^{D_\chi}$. By Requirement~\ref{sec:requireI}, there is finite
  $F$ such that $D \subseteq^* \bigsqcup_{j \in F} R_j \cup
  \bigsqcup_{j \in F} D_j$.  Since $D$ is a $1$-successor of $D_\chi$,
  so is $D-\bigsqcup_{j \in F} R_j$.  Since $D$ and $D_\chi$ are
  disjoint we can assume that if $l(j) = (\chi,0)$ then $j \not\in F$.
  Now if $j \in F$ then $D_j \cap R_{\chi,i} = \emptyset$.
  Contradiction.
\end{proof}

\begin{definition}\mbox{}
  \begin{enumerate}
  \item $D$ is a \emph{$0$-successor} witnessed by $\R^D$ iff $D = A$
    and the lists, $\R^A$ and $\R^D$, are identical.
  \item $D$ is a $1$-successor of $A$ over $\R^A$ was defined in
    Definition~\ref{sec:definable-view-our-6}~(2).
  \item Let $\Dtilde$ be a $n$-successor of $A$ witnessed by
    $\R^W$. If an $\R^{\Dtilde}$ list over $\R^{W}$ exists and $D$ is
    a $1$-successor of $\Dtilde$ over $\R^{\Dtilde}$, then \emph{$D$
      is an $n+1$-successor of $A$ witnessed by $\R^{\Dtilde}$}.
  \item $D$ is a \emph{successor} of $A$ iff, for some $n \geq 0$, $D$
    is an $n$-successor.
  \end{enumerate}
\end{definition}

\begin{corollary}\label{sec:defin-str}
  Let $\chi \in T$.  Then $D_\chi$ is a $|\chi|$-successor of $A$ over
  $\R^{D_{\chi^-}}$.  Furthermore, if $F$ is finite then $D_{\chi} -
  \bigsqcup_{i \in F}R_i$ is a $|\chi|$-successor of $A$ over
  $\R^{D_{\chi^-}}$.
\end{corollary}

\begin{corollary}\label{sec:problemI}
  For all $\chi, \zeta \in T$, if $\chi \neq \zeta$ and $|\chi| =
  |\xi|$, then $D_\chi \sqcup D_\xi$ is a $|\chi|$-successor of $A$
  witnesses by $\R^{D_{\chi^-}}$.
\end{corollary}

We want to transfer these results to the hatted side.  We want to find
$n$-successors of $\Ahat$, \emph{without using the $\Phi$}, witnessing
$A$ and $\Ahat$ are in the same orbit. Just from knowing $A$ and
$\Ahat$ are in the same orbit we want to be able to recover all
successors of $\Ahat$.  But first we need the following lemmas.

\begin{lemma}[Schwarz, see Theorem XII.4.13(ii) of \citet{Soare:87}]
  The index set of maximal sets is $\Pi^0_4$-complete and hence
  computable in $\mathbf{0^{(4)}}$.
\end{lemma}

\begin{lemma}
  The index set of computable sets is $\Sigma^0_3$-complete and hence
  computable in $\mathbf{0^{(3)}}$.
\end{lemma}

\begin{corollary}
  The set $\{\langle e_1, e_2 \rangle : W_{e_1} \text{ lives in }
  W_{e_2}\}$ is $\Sigma^0_5$ and hence computable in $0^{(5)}$.
\end{corollary}

\begin{lemma}\label{sec:definable-view-our-1}
  An $\R^{\Ahat}$ list exists and can be found in an oracle for
  $0^{(5)}$.
\end{lemma}

\begin{proof}
  First we know $\R^{D_\lambda}$ is an $\R^A$ list.  So
  $R^{\Dhat_\lambda}_i = \Phi(R^{D_\lambda}_i)$ is an
  $\R^{\Dhat_\lambda}$ list.  Hence an $\R^{\Ahat}$ list exists.
  However, using $\Phi$ in this fashion does not necessarily bound the
  complexity of $\R^{\Ahat}$.

  Inductively, using an oracle for $0^{(5)}$, we will create an
  $\R^{\Ahat}$ list.  Assume that $\Rhat^{\Ahat}_i$ are known for $i <
  j$, and that for $e < j$, if $\Ahat$ lives in $W_e$ then there is an
  $i < j$ such that $W_e - \Mhat^{W_e} =^* \Rhat^{\Ahat}_i -
  \widehat{M}^{\Rhat^{\Ahat}_i}$.  Look for the least $e \geq j$ such
  that $\Ahat$ lives in $W_e$ and for all $i < j$ such that $W_e -
  \Mhat^{W_e} \neq^* \Rhat^{\Ahat}_i -
  \widehat{M}^{\Rhat^{\Ahat}_i}$. Such an $e$ must exist since an
  $\R^{\Ahat}$ lists exists.  Let $\Rhat^{\Ahat}_j = W_e$.  Apply the
  hatted version of Lemma~\ref{sec:livesI} to get the
  $\Rhat^{\Ahat}_j$ disjoint from $\Rhat^{\Ahat}_i$.
\end{proof}

\begin{definition}
  Let $g$ be such that $W_{g(i)}= R^{\widehat{\Dtilde}}_i$.  Then we
  will say that $g$ is a \emph{presentation} of
  $\R^{\widehat{\Dtilde}}$.
\end{definition}

\begin{lemma}
  Let $\widehat{\Dtilde}$ and an $\R^{\widehat{\Dtilde}}$ list be
  given.  Assume that $g$ is a presentation of
  $\R^{\widehat{\Dtilde}}$. Then all the $1$-successors of
  $\widehat{\Dtilde}$ over $\R^{\widehat{\Dtilde}}$ can be found using
  an oracle for $(g \join 0^{(5)})^{(2)}$.
\end{lemma}

\begin{proof}
  Asking ``whether an $e$ such that $W_e = W_{g(i)}$ and $\Dhat$
  lives in $W_e$'' is computable in $g \join 0^{(5)}$.  $\Dhat$ is a
  $1$-successor of $\widehat{\Dtilde}$ over $\R^{\widehat{\Dtilde}}$
  iff there is a $k$, for all $i \geq k$, [there is an $e$ such that
  $W_e = W_{g(i)}$ and $\Dhat$ lives in $W_e$].
\end{proof}

\begin{corollary}
  The $1$-successors of $\Ahat$ can be found with an oracle for
  $0^{(7)}$.
\end{corollary}

A word of caution: For all $\chi \in T$ of length one, $\Phi(D_\chi)$
is a $1$-successor of $\Ahat$ and, for $\Phi(D_\chi)$, an infinite
$\R^{\Phi(D_\chi)}$ list over $\R^{\Ahat}$ exists.  But, by
Corollary~\ref{sec:problemI}, not every $1$-successor $\Dhat$ of
$\Ahat$ is the image of some such $D_\chi$ even modulo finite many
$R_{\xi,i}$.  Furthermore, there is no reason to believe that if
$\Dhat$ is a $1$-successor of $\Ahat$ that an $\R^{\Dhat}$ list over
$\R^{\Ahat}$ exists.  Unfortunately, we must fix this situation before
continuing.

\begin{definition}
  Let $D_1$ and $D_2$ be $1$-successors of $\Dtilde$ over some
  $\R^{\Dtilde}$ list.  Let an $\R^{D_i}$ list be given. $D_1$ and
  $D_2$ are \emph{$T$-equivalent} iff for almost all $m$, there is an
  $n$ such that $R^{D_1}_m - M^{R^{D_1}_m} =^* R^{D_2}_n -
  M^{R^{D_2}_n}$ and for almost all $m$, there is an $n$ such that
  $R^{D_2}_m - M^{R^{D_2}_m} =^* R^{D_1}_n - M^{R^{D_1}_n}$.
\end{definition}

\begin{lemma}
  If $\chi \in T$ and $F$ is finite then $D_{\chi}$ and $D_{\chi} -
  \bigsqcup_{i \in F}R_i$ are $T$-equivalent $1$-successors of
  $D_{\chi^-}$ over $\R^{D_{\chi^-}}$.
\end{lemma}

\begin{lemma}
  For all $\chi, \zeta \in T$, if $\chi \neq \xi$ and $|\chi| =
  |\xi|$, then $D_\chi$, $D_\xi$ and $D_\chi \sqcup D_\xi$ are
  pairwise $T$-nonequivalent $1$-successors of $D_{\chi^-}$ over
  $\R^{D_{\chi^-}}$.
\end{lemma}

\begin{lemma}
  $D_1$ and $D_2$ are $T$-equivalent iff their automorphic images are
  $T$-equivalent.
\end{lemma}

\begin{lemma}
  Whether ``$\Dhat_1$ and $\Dhat_2$ are $T$-equivalent'' can be
  determined with an oracle for $(g_1 \join g_2 \join \tilde{g} \join
  0^{(5)})^{(2)}$, where $g_i$ and $\tilde{g}$ are representatives of
  needed lists.
\end{lemma}

So $D_\chi$ and $D_\chi-R_i$ are $T$-equivalent. Therefore, we need to
look at the $T$-equivalence class of $D_\chi$ rather than just
$D_\chi$; $D_\chi$ is just a nice representative of the
$T$-equivalence class of $D_\chi$.  $T$-equivalence allows us to
separate $D_\chi$ for $\chi$ of the same length; they are not
$T$-equivalent.  However, we cannot eliminate the image of the
disjoint union of two different $D_\chi$ as a possible successor of
the image of $\Dhat_{\chi^-}$.  For that we need another notion.

\begin{definition}
  Let $D$ be a $1$-successor of $\Dtilde$ over some $\R^{\Dtilde}$
  list.  Let an $\R^{D}$ list be given.  We say that $D$ is
  \emph{atomic} iff for all nontrivial splits $D_1 \sqcup D_2 = D$, if
  $D_i$ is a $1$-successor of $\Dtilde$ then, for almost all $m$,
  $D_i$ lives in $R^{D}_m$.
\end{definition}

\begin{lemma}
  Assume $D$ is an atomic $1$-successor of $\Dtilde$ over some
  $\R^{\Dtilde}$, an $\R^{D}$ list exists, and $D_1 \sqcup D_2$ is a
  nontrivial split of $D$. If $D_i$ is a $1$-successor of $\Dtilde$ then
  an $\R^{D_i}$ list exists and $D$ and $D_i$ are $T$-equivalent.
\end{lemma}

\begin{definition}
  A $T$-equivalent class $\mathcal{C}$ is called an \emph{atomic
    $T$-equivalent class} if every member of $\mathcal{C}$ is atomic.
\end{definition}

The following lemma says that the notion of being atomic indeed
eliminates the disjoint union possibility.

\begin{lemma}
  If $\chi \neq \xi$ and $|\chi| = |\xi|$ then $D_\chi \sqcup D_\xi$
  is not atomic.
\end{lemma}

\begin{lemma}\label{sec:definable-view-our-10}
  Let $D$ be a $1$-successor of $\Dtilde$ over some $\R^{\Dtilde}$
  list.  Let an $\R^{D}$ list be given.  Then $D$ is atomic iff its
  automorphic image is atomic.
\end{lemma}

\begin{lemma}
  Let $\Dhat$ be a $1$-successor of $\widehat{\Dtilde}$ over some
  $\R^{\widehat{\Dtilde}}$ list.  Let an $\R^{\Dhat}$ list be given.
  Determining ``whether $\Dhat$ is atomic'' can be done using an
  oracle for $(g \join \tilde{g} \join 0^{(5)})^{(3)}$, where $g$ and
  $\tilde{g}$ are representatives of needed lists.
\end{lemma}

Unfortunately, with the construction as given so far, there is no
reason to believe that each $D_\chi$ is atomic.  We are going to have
to modify the construction so that each $D_\chi$ is atomic. Thus, we
are going to have to add this as another requirement.

\begin{requirement}\label{sec:requirementI}
  Fix $\chi$ such that $\chi \in T$.  Then $D_\chi$ is an atomic
  $1$-successor of $D_{\chi^-}$.
\end{requirement}

We will have to modify the construction so that we can meet the above
requirement.  This will be done in
Section~\ref{sec:needed-modification}.  Until that section, we will
work under the assumption we have met the above requirement.

These next two lemmas must be proved simultaneously by induction.
They are crucial in that they reduce the apparent complexity down to
something arithmetical.

\begin{lemma}\label{codingstocodings}
  Fix an automorphism $\Phi$ of $\E$ taking $A$ to $\Ahat$. Let
  $\mathcal{C}_{n+1}$ be the class formed by taking all sets of the
  form $\Phi(D_\chi$), where $\chi \in T$ and has length $n+1$, and
  closing under $T$-equivalence.  The collection of all atomic
  $n+1$-successors of $\Ahat$ and $\mathcal{C}_{n+1}$ are the same
  class.
\end{lemma}

\begin{proof}
  For the base case, by Lemma~\ref{sec:definable-view-our-1}, an
  $\R^{\Ahat}$ list exists. Now apply
  Lemma~\ref{sec:definable-view-our-10}. For the inductive case, use
  the following lemma, and then Lemma~\ref{sec:definable-view-our-10}.
\end{proof}

\begin{lemma}\label{sec:definable-view-our-9}
  Let $\widehat{\Dtilde}$ be an atomic $n$-successor of $\Ahat$
  witnessed by $\R^{\What}$. Assume an $\R^{\widehat{\Dtilde}}$ list
  over $\R^{\What}$ exists and $\Dhat$ is an atomic $1$-successor of
  $\widehat{\Dtilde}$ over $\R^{\widehat{\Dtilde}}$. (Then $\Dhat$ is
  an atomic $n+1$-successor of $\Ahat$ witnessed by
  $\R^{\widehat{\Dtilde}}$.)  Then an $\R^{\Dhat}$ list over
  $\R^{\widehat{\Dtilde}}$ can be constructed with an oracle for $g
  \join 0^{(5)}$, where $g$ is representative for
  $\R^{\widehat{\Dtilde}}$.
\end{lemma}

\begin{proof}
  First we will show an $\R^{\Dhat}$ list must exist.  By the above
  lemma, $\Dhat$ is $T$-equivalent to $\Phi(D_\chi)$, where $\chi$ has
  length $n+1$.  An $\R^{D_\chi}$ list exists; hence, so does an
  $\R^{\Dhat}$ list.

  Because of the given properties of $\widehat{\Dtilde}$, the
  $\R^{\What}$ list, and $\R^{\widehat{\Dtilde}}$, if $\Rhat$ is a
  set in the $\R^{\Dhat}$ list, then $\widehat{\Dtilde}$ does not live
  in $\Rhat$.  (This is true for the pre-images of these sets and hence
  for these sets.)
  
  Inductively using an oracle for $g \join 0^{(5)}$ we will create an
  $\R^{\Dhat}$ list.  Assume that $\Rhat^{\Dhat}_i$ are known for $i <
  j$ and that for $e < j$ if $\Dhat$ lives in $W_e$ then there is an
  $i < j$ such that $W_e - \Mhat^{W_e} =^* \Rhat^{\Dhat}_i -
  \widehat{M}^{\Rhat^{\Dhat}_i}$.  Look for the least $e \geq j$ such
  that $\Dhat$ lives in $W_e$, $\widehat{\Dtilde}$ does not live in
  $W_e$, and for all $i < j$ such that $W_e - \Mhat^{W_e} \neq^*
  \Rhat^{\Dhat}_i - \widehat{M}^{\Rhat^{\Dhat}_i}$. Such an $e$ must
  exist. Let $\Rhat^{\Dhat}_j = W_e$.  Apply the hatted version of
  Lemma~\ref{sec:livesI} to get the $\Rhat^{\Ahat}_j$ disjoint from
  $\Rhat^{\Ahat}_i$.
\end{proof}

\begin{definition}\label{sec:TA}
  Let $\mathcal{T}(A)$ denote the class of atomic $T$-equivalence
  classes of successors (of $A$) with the binary relation of
  $1$-successor restricted to successors of $A$.
\end{definition}

\begin{remark}\label{sec:LambdaI}
  So the map $\Lambda(\chi) = D_\chi$ is a map from $T$ to
  $\mathcal{T}(A)$ taking a node to a representative of an atomic
  $T$-equivalent class of successors.  Furthermore, $\zeta$ is an
  immediate successor of $\chi$ iff $D_\zeta$ is a $1$-successor of
  $D_\chi$.  Hence $\Lambda$ is an isomorphism. Recall $D_\chi =
  D_\alpha$ if $l(\alpha) = (\chi,0)$.  Hence $\Lambda$ is computable
  along the true path which is computable in $\mathbf{0^{(2)}}$.
\end{remark}

\begin{lemma}\label{gameplan2}
  If $A$ and $\Ahat$ are in the same orbit witnessed by $\Phi$ then
  $\mathcal{T}(\Ahat)$ must exist and must be isomorphic to
  $\mathcal{T}(A)$ via an isomorphism induced by $\Phi$ and computable
  in $\Phi \join 0''$.  The composition of this induced isomorphism
  and the above $\Lambda$ is an isomorphism between
  $\mathcal{T}(\Ahat)$ and $T$.  (This addresses part two of our game
  plan.)
\end{lemma}

Our coding is not elementary; it is not even in
$\mathcal{L}_{\omega_1,\omega}$.  The coding depends on the infinite
lists $\R^{D\chi}$.  One cannot say such a list exists in
$\mathcal{L}_{\omega_1,\omega}$.  It is open if there another coding
of $T$ which is elementary or in $\mathcal{L}_{\omega_1,\omega}$.
This is another excellent open question.

\begin{lemma}\label{oracle}
  $\mathcal{T}(\Ahat)$ has a presentation computable in
  $\mathbf{0^{(8)}}$.
\end{lemma}

\subsection{More Requirements; The Homogeneity Requirements}
\label{manytrees}

Let $\chi, \xi \in T$ be such that $\chi^- = \xi^-$.  Then in terms of
the above coding, the atomic $T$-equivalence classes of $D_\chi$ and
$D_\xi$ cannot be differentiated.  For almost all $i$, $D_\chi$ and
$D_\xi$ live in $R_{\chi^-,2i}$ (if $M_{\chi^-,i}$ is maximal in
$R_{\chi^-,i}$) and
\begin{equation}
  \label{eq:9}
  \text{for all } i (D_\chi \text{ lives in  } 
  R_{\chi,i} \text{ iff } D_\xi \text{ lives in } R_{\xi,i}).
\end{equation}
In this sense, these sets are homogeneous.  What we are about to do
has the potential to destroy this homogeneity.  We must be careful not
to destroy this homogeneity.

In fact, we must do far more than just restore this homogeneity. For
each $T_i$ we will construct an $A_{T_i}$.  For all $\chi^{T_k} \in
T_k$, we will construct $D_{\chi^{T_k}}$, $R_{\chi^{T_k},i}$, and
$M_{\chi^{T_k},i}$. In order to complete part~3 of our game plan (that
is, sets coded by isomorphic trees belong to the same orbit) we must
ensure that the following homogeneity requirement holds.

\begin{requirement}\label{sec:homo}
  For all $k, \widehat{k}$, if $\chi^{T_{k}} \in T_k$,
  $\chi^{T_{\widehat{k}}} \in T_{\widehat{k}}$, and $|\chi^{T_k}| = |
  \chi^{T_{\widehat{k}}}|$ then, for all $i$,
  \begin{equation*}\label{eq:homo1}
    M_{\chi^{T_k},i} \text{ is maximal in }  R_{\chi^{T_k},i} \text{ iff }
    M_{\chi^{T_{\widehat{k}}},i} \text{ is maximal in } 
    R_{\chi^{T_{\widehat{k}}},i}, \text{ and}
  \end{equation*}
  \begin{equation*}\label{eq:homo2}
    M_{\chi^{T_k},i} =^*  R_{\chi^{T_k},i}  \text{ iff } 
    M_{\chi^{T_{\widehat{k}}},i} =^*  R_{\chi^{T_{\widehat{k}}},i}. 
  \end{equation*}
\end{requirement}

\begin{remark}
  We cannot overstate the importance of this requirement. It is key to
  the construction of \emph{all} of the needed automorphisms; see
  Section~\ref{sec:defining--phi}. Note that we use
  Section~\ref{sec:defining--phi} twice; once in this proof and once
  in the proof of Theorem~\ref{sec:maincor}. 
\end{remark}

One consequence of this requirement is that we must construct all the
sets, $D_{\chi^{T_k}}$, $R_{\chi^{T_k},i}$, and $M_{\chi^{T_k},i}$,
simultaneously using the same tree of strategies.  Up to this point we
have been working with a single $T$.  To dovetail all the trees into
our construction at the node $\alpha \in Tr$ where $|\alpha| = k$ we
will start coding tree $T_k$. Since at each node we only needed
answers to a finite number of $\Delta^0_3$ questions, this dovetailing
is legal in terms of the tree argument.  Note that each tree $T$ gets
its own copy of $\omega$ to work with.

So at each $\alpha \in Tr$, we will construct, for $k < |\alpha|$,
$R_\alpha^k, E^k_\alpha, M^k_\alpha$, and $D_\alpha^k$ as above.  The
$e$th marker for $M^k_\alpha$ be will denoted
$\Gamma^{\alpha,k}_e$. Assume that $\alpha \subset f$, $k < |\alpha|$,
and $l(|\alpha|-k) = (\chi,i)$; then $R^k_\alpha = R_{\chi^{T_k},2i}$,
$M^k_\alpha = M_{\chi^{T_k},2i}$, $\Gamma^{\alpha,k}_e =
\Gamma^{\chi^{T_k},2i}_e$, $E_\alpha = R_{\chi^{T_k},2i+1} =
M_{\chi^{T_k},2i+1}$, and if $i = 0$ then $D_\alpha = D_{\chi^{T_k}}$.
In the following, when the meaning is clear, we drop the subscript $k$
and assume we are working with a tree $T$.

\subsection{Meeting the Remaining
  Requirements}\label{sec:needed-modification}

The goal in this section is to understand what it takes to show
$D_\chi$ is atomic when $\chi \in T \cup \{\lambda\}$.  We have to do
this and meet Requirements~\ref{require3} and \ref{sec:homo}.  Since
we will meet Requirement~\ref{require3}, an $\R^{D_\chi}$ list exists.
$D_\chi$ is potentially not atomic witnessed by an c.e.\ set $W
\subseteq D_\chi$ if the set of $i$ such that $W$ lives in an
$R^{D_\chi}_i$ is an infinite coinfinite set.  We must make sure that
$W$ behaves cohesively on the sets $R^{D_\chi}_i$.

We will meet Requirement~\ref{sec:requirementI} by a $e$-state
argument on the $R_{\chi,2i}$s; this is similar to a maximal set
construction.  With the maximal set construction, for each $e, s$
and $x$, there are $2$ states, either state $0$ iff $x \not\in W_{e,s}$
or $1$ iff $x \in W_{e,s}$.  Here the situation is more complex.

$R$ has state $1$ w.r.t.\ a single $e$ iff $W_e \cap R =^* \emptyset$.
$R$ has state $2$ w.r.t.\ a single $e$ iff $W_e \cap R \neq^*
\emptyset$.  $R$ has state $3$ w.r.t.\ a single $e$ iff there is an
$\breve{e}$ such that $W_e \sqcup W_{\breve{e}} = M^R$. $R$ has state
$4$ w.r.t.\ a single $e$ iff there is an $\breve{e}$ such that $W_e
\sqcup W_{\breve{e}} = R$.  If the highest state of $R$ is $3$ w.r.t.\
a single $e$ then $W_e$ is a nontrivial split of $M^R$. Determining
the state of $R$ w.r.t.\ $e$ is $\Sigma^0_3$.  (The state $0$ will be
used later.)

Let $s_{e'}$ be the state of $R$ w.r.t.\ $e'$. The $e$-state of $R$ is
the string $s_0s_1s_2 \ldots s_e$.  An $e$-state $\sigma_1$ is
\emph{greater} than $\sigma_2$ iff $\sigma_1 <_L \sigma_2$.  We will
do an $e$-state construction along the true path for the tree $T_k$.

Assume $\alpha \in Tr$, $e = |\alpha| -k$, and $l(e) = (\chi,i)$.
Since at $\alpha$ we can get answers to a finite number of
$\Delta^0_3$ questions, at $\alpha$ we will have encoded answers to
which, if any, $\beta \prec \alpha$, if $l(|\beta|-k) = (\chi,i')$ is
$W_e \cap R^k_{\beta}$ is infinite; for which of the above $\beta$s
and for which $j < e$, does $W_l$, for $l < e$, witness that $W_j$ is
a split of $M^k_{\beta}$;  and for which of the above $\beta$s, and
for which $j<e$, does $W_l$, for $l < e$, witness that $W_j$ is a
split of $R^k_{\beta}$?

Using this information, $\alpha$ will determine $\beta^{\alpha,k}_0,
\beta^{\alpha,k}_1, \ldots$ such that $R^k_{\beta^{\alpha,k}_e}$ has
the greatest possible $e$-state according to the information encoded at
$\alpha$.  This listing does not change w.r.t.\ stage.  For all other
$\beta \prec \alpha$ such that $l(|\beta|-k) = (\chi,i')$, when
$\alpha \subseteq f_s$, $\alpha$ will dump
$\Gamma^{\beta,k}_{|\alpha|}$ into $M^k_{\beta}$.  If $\alpha \subset
f$ then, for the above $\beta$, $R^k_{\beta} =^* M^k_{\beta}$.

One can show, for each $e$, there is an $\alpha_e \subset f$ such that,
for all $\gamma$ with $\alpha_e \preceq \gamma \subset f$,
$\beta^{\alpha_e,k}_e =\beta^{\gamma,k}_e $. Hence $M^k_{\alpha_e}$ is
maximal in $R^k_{\alpha_e}$ and Requirement~\ref{require3} is met.
In addition, one can show that, for almost all $i$, such that
$M_{\chi^{T_k},i} \neq ^* R_{\chi^{T_k},i}$, $R_{\chi^{T_k},i}$ have
the same $e$-state and hence $D_\chi$ is atomic.

However, Equation~\ref{eq:9} no longer holds and hence
Requirement~\ref{sec:homo} is not met. The problem is that we can
dump $M_{\chi^{T_k},i}$ into $R_{\chi^{T_k},i}$ without dumping
$M_{\chi^{T_{\widehat{k}}},i}$ into $R_{\chi^{T_{\widehat{k}}},i}$.

The solution is that when we dump $M_{\chi^{T_k},i}$ we also must dump
$M_{\chi^{T_{\widehat{k}}},i}$, for all possible
$\chi^{T_{\widehat{k}}}$.  This means that we have to do the above
$e$-state construction for $T_k$ simultaneously for all $T_k$.  So,
for each $n$, we have \emph{one} $e$-state construction, for all
$D_{\chi^{T_k}}$ and $D_{\xi^{T_k}}$, for all $k$ and for all
$\chi^{T_k}, \xi^{T_k} \in T_k$ with $|\chi^{T_k}| = |\xi^{T_k}| =n$.

To do this we need the following notation: Let $\{\xi_i: i \in
\omega\}$ be a computable listing of all nodes of length $n$ in
$\omega^{< \omega}$.  Fix some nice one-to-one onto computable
listing, $ \langle -, -,-\rangle$, of all triples $( e, k, l )$ and,
furthermore, assume if $( e, k, l )$ is the $m$th triple listed then
$\langle e, k, l \rangle = m$.

Assume $l(|\beta|) = (\xi,i)$ and $|\xi| = n$. If there is a $\beta'
\preceq \beta$ and a $k \leq |\beta|$ such that $l(|\beta'|-k) =
(\xi',i)$ (the same $i$ as above) and $|\xi'| = n$ and, furthermore,
$\beta'$ is the $l$th such $\beta'$ then the state of $\beta$
w.r.t. $\langle e, k, l \rangle$ is the state of $R^k_{\beta'}$ w.r.t.\
$e$.  Otherwise the state of $\beta$ w.r.t.\ $\langle e, k, l \rangle$
is $0$.  Let $s_{\langle e', k', l' \rangle}$ be the state of $\beta$
w.r.t.\ $\langle e', k', l' \rangle$. The $\langle e, k, l
\rangle$-state of $\beta$ is the string $s_{\langle e_0, k_0, l_0
  \rangle}s_{\langle e_1, k_1, l_1 \rangle}s_{\langle e_2, k_2, l_2
  \rangle} \ldots s_{\langle e, k, l \rangle}$.

Using the additional information we encoded into $\alpha$ for the
single $e$-state construction, $\alpha$ has enough information to
determine the $\langle e, k, l \rangle$-state of $\beta \preceq
\alpha$.  Using this information, $\alpha$ will determine
$\beta^{\alpha}_{\langle e_0, k_0, l_0 \rangle},
\beta^{\alpha}_{\langle e_1, k_1, l_1 \rangle}, \ldots$ such that
${\beta^{\alpha}_{\langle e, k, l \rangle}}$ has the greatest possible
$\langle e, k, l \rangle$-state according the information encoded at
$\alpha$.  Again this listing does not change w.r.t.\ stage.

For all other $\beta \prec \alpha$ such that $l(|\beta|) =
(\xi_{j'},i')$ and $|\xi_{j'}| = n$, when $\alpha \subseteq f_s$, for
all $k$, for all $\chi \in T^k$ of length $n$, $\alpha$ will dump
$\Gamma^{\chi^{T_k},2i'}_{|\alpha|}$ into $M_{\chi^{T_k},2i'}$.  If
$\alpha \subset f$ then, for the above $i'$, for all $k$, for all
$\chi \in T^k$ of length $n$, $M_{\chi^{T_k},2i'} =^*
R_{\chi^{T_k},2i'}$.

One can show, for each $\langle e, k, l \rangle$, there is an
$\alpha_{\langle e, k, l \rangle} \subset f$ such that for all
$\gamma$ with $\alpha_{\langle e, k, l \rangle} \preceq \gamma \subset
f$, $\beta^{\alpha_{\langle e, k, l \rangle}}_{\langle e, k, l
  \rangle} =\beta^{\gamma}_{\langle e, k, l \rangle}$.  Assume
$l(|\beta^{\alpha_{\langle e, k, l \rangle}}_{\langle e, k, l
  \rangle}|) = ( \chi,i)$.  Then, for all $k$, for all $\chi^{T_k} \in
T^k$ of length $n$, $M_{\chi^{T_k},2i}$ is maximal in
$R_{\chi^{T_k},2i}$.  Hence Requirements~\ref{require3} and
~\ref{sec:homo} are met.

In addition, one can show that, for all $e$, for all $k$, for all
$\chi^{T_k} \in T^k$ of length $n$, for almost all $i$, if
$M_{\chi^{T_k},2i'}$ is maximal in $R_{\chi^{T_k},2i'}$ then
$R_{\chi^{T_k},i}$ has the same state w.r.t.\ $e$ and, hence,
$D_{\chi^{T_k}}$ is atomic.  Thus Requirement~\ref{sec:requirementI}
is met.

\subsection{Same Orbit}\label{sec:same-orbit}

Let $T$ and $\widehat{T}$ be isomorphic trees via an isomorphism
$\Lambda$. We must build an automorphism $\Phi_\Lambda$ of $\E$ taking
$A$ to $\Ahat$.  We want to do this piecewise.  That is, we want to
build isomorphisms between the $\E^*(D_\chi)$ and
$\E^*(\Dhat_{\Lambda(\chi)})$ and piece them together in some fashion
to get an automorphism.  Examples of automorphisms constructed in such
a manner can be found in Section~5 of \citet{mr2001k:03086} and
Section~7 of \citet{Cholak.Harrington:nd}.

In reality $T=T_k$ and $\That = T_{\widehat{k}}$. The sets in question
for $T_k$ are $D_{\chi^{T_k}}$, $R_{\chi^{T_k},i}$, and
$M_{\chi^{T_k},i}$.  Here we will just drop the $T_{\widehat{k}}$
superscript from $\chi$.  The sets in question for $T_{\widehat{k}}$
are $D_{\chi^{T_{\widehat{k}}}}$, $R_{\chi^{T_{\widehat{k}}},i}$, and
$M_{\chi^{T_{\widehat{k}}},i}$. Here we will ``hat'' the sets involved
and drop the $T_{\widehat{k}}$ superscript  from $\chi$.

However, before we shift to our standard notation changes we would
like to point out  the following. Since $\Lambda$ is an isomorphism between
$T_k$ and $T_{\widehat{k}}$, $|\chi^T_k| = |\Lambda(\chi^{T_k})|$.
Therefore, by Requirement~\ref{sec:homo}, for all $i$,
\begin{equation*}\label{eq:homo3}
  M_{\chi^{T_k},i} \text{ is maximal in }  R_{\chi^{T_k},i} \text{ iff }
  M_{\Lambda(\chi^{T_k}),i} \text{ is maximal in } 
  R_{\Lambda(\chi^{T_k}),i},
\end{equation*}
\begin{equation*}\label{eq:homo4}
  \text{and } M_{\chi^{T_k},i} =^*  R_{\chi^{T_k},i}  \text{ iff } 
  M_{\Lambda(\chi^{T_k}),i} =^*  R_{\Lambda(\chi^{T_k}),i}. 
\end{equation*}

\subsubsection{Extendible algebras of computable
  sets}\label{sec:oldwork}

The workhorse for constructing $\Phi_\Lambda$ is the following theorem
and two lemmas.

\begin{theorem}[Theorem~5.10 of
  \citet{Cholak.Harrington:nd}]\label{interface}
  Let $\B$ be an extendible algebra of computable sets and similarly
  for $\widehat{\B}$.  Assume the two are extendibly isomorphic via
  $\Pi$.  Then there is a $\Phi$ such that $\Phi$ is a $\Delta^0_3$
  isomorphism between $\E^*(A)$ and $\E^*(\widehat{A})$, $\Phi$ maps
  computable subsets to computable subsets, and, for all $R \in \B$,
  $(\Pi(R) - \widehat{A}) \sqcup \Phi(R \cap A)$ is computable\ (and
  dually).
\end{theorem}

\begin{lemma}\label{niceextendiblealgebra}
  Let $\chi \in T$. The collection of all $R_{\chi,i}$ forms an
  extendible algebra, $\mathcal{B}_\chi$, of computable sets.
\end{lemma}

\begin{proof}
  Apply Theorem~2.17 of \citet{Cholak.Harrington:nd} to $A = \omega$
  to get an extendible algebra of $\S_\R(\omega)$ of all computable
  sets with representation $B$.  Let $j \in B_\chi$ iff there is an $i
  \leq j$ such that $S_j = R_{\chi,i}$.  Now take the subalgebra
  generated by $B_\chi$ to get $\mathcal{B}_\chi$.
\end{proof}

\begin{lemma}
  Let $\chi \in T$; then the join of $\B_{\chi^-}$ and $\B_\chi$ is an
  extendible algebra of computable sets, $\B_{\chi^- \join \chi}$.
\end{lemma}

\begin{proof}
  See Lemma 2.16 of \citet{Cholak.Harrington:nd}.
\end{proof}

\begin{lemma}
  For all $i$, if $R_{\xi,j} \not\equiv_\R R_{\chi,i}$ and $R_{\xi,j}
  \not\equiv_\R R_{\chi^-,i}$ then $D_\chi \cap R_{\xi,j} =
  \emptyset$.
\end{lemma}

\begin{proof}
  See Lemma~\ref{sec:definable-view-our-5}.
\end{proof}

\begin{lemma}\label{isomorphismnice}
  If $\chi,\xi \in T$ and $|\chi| = |\xi|$ then $\B_{\chi^- \join
    \chi}$ and $\widehat{\B}_{\xi^- \join \xi}$ are extendibly
  isomorphic via $\Phi_{\chi,\xi}$ where
  $\Phi_{\chi,\xi}(R_{\chi^-,i}) =\Rhat_{\xi^-,i}$ and
  $\Phi_{\chi,\xi}(R_{\chi,i}) =\Rhat_{\xi,i}$.  Furthermore,
  $\Phi_{\chi,\xi}$ is $\Delta^0_3$.
\end{lemma}

\subsubsection{Building $\Phi_\Lambda$ on the $D$s and $M$s}
\label{sec:build-phiI}

The idea is to use Theorem~\ref{interface} to map $\E^*(D_\chi)$ to
$\E^*(\Dhat_{\Lambda(\chi)})$.  By the above lemmas, there is little
question that the extendible algebras we need are some nice
subalgebras of $\B_{\chi^- \join \chi}$ and
$\widehat{\B}_{\Lambda(\chi^-) \join \Lambda(\chi)}$ and the
isomorphism between these nice subalgebras is induced by the
isomorphism $\Phi_{\chi,\Lambda(\chi)}$.

We will use the following stepwise procedure to define part of
$\Phi_\Lambda$.  This is not a computable procedure but computable in
$\Lambda \join 0''$. $\chi$ is added to $\N$ at step $s$ iff we
determined the image of $D_\chi$ (modulo finitely many
$R_{\chi^-,j}$). The parameter $i_{\chi,s}$ will be used to keep track
of the $M_{\chi,i}$ which we have handled and will be increasing
stepwise.  This procedure does not completely define $\Phi_\Lambda$;
we will have to deal with those $W$ which are not subsets of
$\bigsqcup M \cup \bigsqcup D$.

\emph{Step $0$}: Let $\mathcal{N}_0 = \{\lambda\}$. By the above
lemmas $\B_\lambda$ is isomorphic to $\widehat{\B}_\lambda$ via
$\Phi_{\lambda,\lambda}$.  Let $i_{\lambda,0} = 0$. Now apply
Theorem~\ref{interface} to define $\Phi_\Lambda$ for $W \subseteq A =
D_\lambda$ and dually.

\emph{Step $s+1$}: \emph{Part $\chi \cat s$}: For each $\chi \in
\N_{s}$ such that $\chi \cat s \in T$ do the following: Add $\chi \cat
s$ to $\N_{s+1}$.  Let $i_{\chi \cat s,s+1} = 0$.  Apply
Lemma~\ref{sec:definable-view-our} to $\chi \cat s$ to get $i'$. Apply the
hatted version of Lemma~\ref{sec:definable-view-our} to $\Lambda(\chi
\cat s)$ to get $\widehat{i}'$.  Let $i_{\chi,s+1}$ be the max of
$i'$, $\widehat{i}'$ and $i_{\chi,s}+1$.  Let $\B^*_{\chi, \chi \cat
  s}$ be the extendible algebra generated by $R_{\chi,i}$, for $i \geq
i_{\chi,s+1}$, and, for all $j$, $R_{\chi \cat s, j}$.  Define
$\B^*_{\Lambda(\chi), \Lambda(\chi \cat s)}$ in a dual fashion.  Now
$\Phi_{\chi \cat s, \Lambda(\chi \cat s)}$ induces an isomorphism
between these two extendible algebras.  Now apply
Theorem~\ref{interface} to define $\Phi_\Lambda$ for $W \subseteq
\big ( D_{\chi \cat s} - \bigsqcup_{i < i_{\chi,s+1}} R_{\chi,i} \big
)$ and $\Phi_\Lambda^{-1}$ for $\What \subseteq \big (
\Dhat_{\Lambda(\chi \cat s)} - \bigsqcup_{i < i_{\chi,s+1}}
\Rhat_{\Lambda(\chi),i} \big ) $.

\emph{Step $s+1$}: \emph{Part $i_{\chi,s+1}$}: For all $\chi \in \N_s$
and for all $i$ such that $i_{\chi,s} \leq i < i_{\chi,s+1}$, do the
following: Let $S_{\chi,i} = \big ( M_{\chi,i} - \bigsqcup_{\xi \in
  \N_s} D_\xi \big )$ and $\Shat_{\Lambda(\chi),i} = \big (
\Mhat_{\chi,i} - \bigsqcup_{\xi \in \N_s} \Dhat_{\Lambda(\xi)} \big
)$. So $H_{\chi,i} \subseteq S_{\chi,i}$ and $\Hhat_{\chi,i} \subseteq
\Shat_{\Lambda(\chi),i}$.  $S_{\chi,i}$ and $\Shat_{\Lambda(\chi),i}$
are both infinite and furthermore, by Equation~\ref{eq:9}, the one is
computable iff the other is computable.

\emph{Subpart $H$}: If both $S_{\chi,i}$ and $\Shat_{\Lambda(\chi),i}$
are noncomputable then apply Theorem~\ref{interface} (using the empty
extendible algebras) to define $\Phi_\Lambda$ for $W \subseteq
S_{\chi,i}$ and $\Phi_\Lambda^{-1}$ for $\What \subseteq
\Shat_{\Lambda(\chi),i}$.  If both $S_{\chi,i}$ and
$\Shat_{\Lambda(\chi),i}$ are computable then such $\Phi_\Lambda$ can
be found by far easier means.

One can show that $T = \lim_s \N_s$ and that, for all $i$, $\chi \in
T$, there is step $s$ such that, for all $t \geq s$, $i_{\chi,t} \geq
i$.  For all $\chi \in T$, let $s_\chi$ be the step that $\chi$ enters
$\N$ and $s_{\chi,i}$ be the first stage such that
$i_{\chi,s_{\chi,i}} > i$.

\subsubsection{Defining $ \Phi_\Lambda$ on $R_{\chi,i}$}
\label{sec:defining--phi}

Let $s = s_{\chi,i}$. By Section~\ref{sec:build-phiI}, $\Phi_\Lambda$
is defined on \[M_{\chi,i} = S_{\chi,i} \sqcup \bigsqcup_{\xi \in
  \N_s} (R_{\chi,i} \cap D_\xi); \]
\[\Phi_\Lambda(M_{\chi,i}) = \Shat_{\Lambda(\chi),i} \sqcup
\bigsqcup_{\xi \in \N_s} \Phi_\Lambda(R_{\chi,i} \cap D_\xi).\] Hence
$\Phi_\Lambda$ is defined on subsets $W$ of $M_{\chi,i}$. Furthermore,
if such a $W$ is computable so is $\Phi_\Lambda(W)$.

Let $\xi \in \N_s$.  Then \[R_{\chi,i} \cap \bigg ( D_\xi -
\bigsqcup_{j < i_{\xi^-,s_\xi}} R_{\xi^-,j}\bigg ) = R_{\chi,i} \cap
D_\xi,\] \[\Rhat_{\Lambda(\chi),i} - \bigg ( \Dhat_{\Lambda(\xi)} -
\bigsqcup_{j < i_{\xi^-,s_\xi}} \Rhat_{\Lambda(\xi^-),j}\bigg) =
\Rhat_{\Lambda(\chi),i} - \Dhat_{\Lambda(\xi)},\] and
$\Phi_{\chi,\Lambda(\chi)}(R_{\chi,i}) = \Rhat_{\Lambda(\chi),i}$.
Therefore, by Theorem~\ref{interface},
\[\Rhat_{\Lambda(\chi),i} - \Dhat_{\Lambda(\xi))} \sqcup
\Phi_\Lambda(R_{\chi,i} \cap D_\xi) = \widehat{X}_\xi\] is computable.
Since $\Phi_\Lambda(R_{\chi,i} \cap D_\xi) \subset
\Dhat_{\Lambda(\xi)}$, $\Rhat_{\Lambda(\chi),i} \triangle
\widehat{X}_\xi \subseteq \Dhat_{\Lambda(\xi)}$ (recall $\triangle$ is
the symmetric difference between two sets).  Fix computable sets
$\Rtilde^{in}_\xi$ and $\Rtilde^{out}_\xi$ such that $\widehat{X}_\xi
= \big ( \Rhat_{\Lambda(\chi),i} \sqcup \Rtilde^{in}_\xi \big ) -
\Rtilde^{out}_\xi$.

Consider the computable set
\[ 
\Rtilde = \big ( \Rhat_{\Lambda(\chi),i} \sqcup \bigsqcup_{\xi \in
  \N_s} \Rtilde^{in}_\xi \big ) - \bigsqcup_{\xi \in \N_s}
\Rtilde^{out}_\xi.\] Then \[ \Rtilde - \bigsqcup_{\xi \in \N_s}
\Phi_\Lambda(R_{\chi,i} \cap D_\xi) = \Shat_{\Lambda(\chi),i} \sqcup
\bigg ( \Rhat_{\Lambda(\chi),i} -M_{\Lambda(\chi),i} \bigg ). 
\]
Therefore
\[ 
\Rtilde - \Phi_\Lambda(M_{\chi,i}) = \Rhat_{\Lambda(\chi),i}
-M_{\Lambda(\chi),i}. 
\]

Since $M_{\chi,i}$ is maximal in $R_{\chi,i}$ or $M_{\chi,i} =^*
R_{\chi,i}$, if $W \subseteq R_{\chi,i}$ either $W \subseteq^*
M_{\chi,i}$ or there is computable $R$ such that $R \subseteq
M_{\chi,i}$ and $R_W \cup W = R_{\chi,i}$.  In the former case,
$\Phi_\Lambda(W)$ is defined.  In the latter case, let
\[ \Phi_\Lambda(W) = (\Rtilde - \Phi_\Lambda(R_W)) \sqcup
\Phi_\Lambda(W \cap R_W).\] Hence $\Phi_\Lambda(R_{\chi,i}) =
\Rtilde$.

Since $\Lambda$ is an isomorphism between $T$ and $\That$, $|\chi| =
|\Lambda(\chi)|$.  Therefore, as we noted above, by
Requirement~\ref{sec:homo}, either $M_{\chi,i}$ is maximal in
$R_{\chi,i}$ and $\Mhat_{\chi,i}$ is maximal in
$\Rhat_{\Lambda(\chi),i}$ or $M_{\chi,i}=^*R_{\chi,i}$ and
$\Mhat_{\chi,i}=^*\Rhat_{\Lambda(\chi),i}$.  In either case,
$\Phi_\Lambda$ induces an isomorphism between $\E^*(R_{\chi,i})$ and
$\E^*(\Rtilde)$.  $\Phi^{-1}_\Lambda$ on
$\E^*(\Rhat_{\Lambda(\chi),i})$ is handled in the dual fashion.

\subsubsection{Putting $ \Phi_\Lambda$ together}

By Requirement~\ref{sec:requireI} and our construction, for all $e$,
there are finite sets $F_D$ and $F_R$ such that either
\begin{equation}
  \label{eq:14}
  W_e \subseteq^* \biggl ( \bigsqcup_{\chi \in F_D}D_\chi \cup
  \bigsqcup_{(\chi,i) \in F_R} R_{\chi,i} \biggr )
\end{equation}
or there is an $R_{W_e}$ such that
\begin{equation}
  \label{eq:7}
  R_{W_e} \subseteq \bigg ( \bigsqcup_{\chi \in F_D}D_\chi \cup
  \bigsqcup_{(\chi,i) \in F_R} R_{\chi,i} \bigg) \text{ and } W_e \cup
  R_{W_e} = \omega.
\end{equation}
It is possible to rewrite the set
\begin{equation*}
  \bigsqcup_{\chi \in F_D}D_\chi \cup
  \bigsqcup_{(\chi,i) \in F_R} R_{\chi,i} 
\end{equation*}
as
\begin{equation}\label{eq:10}
  \bigsqcup_{\chi \in F_D} \bigg ( D_\chi - \bigsqcup_{(\xi,j) \in
    F_\chi} R_{\xi,j} \bigg ) \sqcup \bigsqcup_{(\chi,i) \in F^*_R}
  R_{\chi,i},
\end{equation}
where $F^*_R \subseteq F_R \cup \bigcup_{\chi \in F_D} F_\chi$ and
$F_\chi$ is finite and includes the set $\{ (\chi^-, l): l <
i_{\chi^-,s_\chi}\}$.  $\Phi_\Lambda$ as defined in the
Section~\ref{sec:build-phiI} is well behaved on the first union in
Equation~\ref{eq:10} and, furthermore, on these unions computable sets
are sent to computable sets.  Similarly, by
Section~\ref{sec:defining--phi}, $\Phi_\Lambda$ is well behaved on the
second union in Equation~\ref{eq:10} and, furthermore, on these unions
computable sets are sent to computable sets.

If Equation~\ref{eq:14} for $e$ hold, then $\Phi(W_e)$ is
determined. Otherwise Equation~\ref{eq:7} holds and map $W_e =
\overline{R_{W_e}} \sqcup (W \cap R_{W_e})$ to
$\overline{\Phi(R_{W_e})} \sqcup \Phi(W \cap R_{W_e})$.
$\Phi^{-1}_\Lambda$ is handled in the dual fashion. So $\Phi_\Lambda$
is an automorphism.

\section{Invariants and Properly $\Delta^0_\alpha$ orbits}
\label{sec:last}

It might appear that $\mathcal{T}(A)$ is an invariant which determines
the orbit of $A$. But there is no reason to believe for an arbitrary
$A$ that $\mathcal{T}(A)$ is well defined.  The following theorem
shows that $\mathcal{T}(\Ahat)$ is an invariant as far as the orbits
of the $A_T$s are concerned.  In Section~\ref{tinv}, we prove a more
technical version of the following theorem.

\begin{theorem}\label{sec:invariant}
  If $\Ahat$ and $A_T$ are automorphic via $\Psi$ and $T \cong
  \mathcal{T}(\Ahat)$ via $\Lambda$ then $A_T \approx \Ahat$ via
  $\Phi_\Lambda$ where $\Phi_\Lambda\leq_T \Lambda \join
  \bf{0}^{(8)}$.
\end{theorem}

\begin{proof}
  See Section~\ref{sec:invarproof}.
\end{proof}

\begin{theorem}[Folklore\footnote{See Section~\ref{folkloreproofII} for
    more information and a proof.}]\label{sec:folkloredark}
  For all 
  finite $\alpha$ there is a computable tree $T_{i_\alpha}$ from the
  list in Theorem~\ref{folkloreI} such that, for all computable trees
  $T$, $T$ and $T_{i_\alpha}$ are isomorphic iff $T$ and
  $T_{i_\alpha}$ are isomorphic via an isomorphism computable in
  $\text{deg}(T) \join 0^{(\alpha)}$.  But, for all $\beta < \alpha$
  there is an $i^*_\beta$ such that $T_{i^*_\beta}$ and $T_{i_\alpha}$
  are isomorphic but are not isomorphic via an isomorphism computable
  in $0^{(\beta)}$.
\end{theorem}

It is open if the above theorem holds for all $\alpha$ such that
$\omega \geq \alpha < \wock$.  But if it does then so does the theorem
below.

\begin{theorem}
  For all finite $\alpha > 8$ there is a properly $\Delta^0_{\alpha}$
  orbit.
\end{theorem}

\begin{proof}
  Assume that $A_{T_{i_\alpha}}$ and $\Ahat$ are automorphic via an
  automorphism $\Phi$.  Hence, by part 2 of the game plan,
  $\mathcal{T}(\Ahat)$ and $T_{i_\alpha}$ are isomorphic.  Since
  $\mathcal{T}(\Ahat)$ is computable in $0^{(8)}$, $\alpha > 8$, and
  by Theorem~\ref{sec:folkloredark}, $\mathcal{T}(\Ahat)$ and
  $T_{i_\alpha}$ via a $\Lambda \leq_T 0^{(\alpha)}$.  By
  Theorem~\ref{sec:invariant}, $\Ahat$ and $A_{T_{i_\alpha}}$ are
  automorphic via an automorphism computable in $0^{(\alpha)}$.

  Fix $\beta$ such that $8 \geq \beta < \alpha$.  By part 3 of the
  game plan and the above paragraph, $A_{T_{i_\alpha}}$ and
  $A_{T_{i^*_\beta}}$ are automorphic via an automorphism computable
  in $0^{(\alpha)}$.  Now assume $A_{T_{i^*_\beta}} \approx
  A_{T_{i_\alpha}}$ via $\Phi$.  By Lemma~\ref{gameplan2},
  $\mathcal{T}(A_{T_{i^*_\beta}}) \cong T_{i_\alpha}$ via
  $\Lambda_{\Phi}$, where $\Lambda_\Phi \leq_T \Phi \join
  \bf{0}^{(2)}$.  Since $\mathcal{T}(A_{T_{i^*_\beta}})$ is computable
  in $0^{(8)}$ and $\mathcal{T}(A_{T_{i^*_\beta}})$ is isomorphic to
  $T_{i^*_\beta}$ via an isomorphism computable in $0^{(\beta)}$ (part
  1 of the game plan), by Theorem~\ref{sec:folkloredark},
  $\Lambda_\Phi >_T 0^{(\beta)}$.  Hence $\Phi >_T 0^{(\beta)}$.
\end{proof}

\subsection{Proof of
  Theorem~\ref{sec:invariant}}\label{sec:invarproof}

For $A_T$ the above construction gives us a $\mathbf{0''}$ listing of
the sets $D_\chi$, $R_{\chi,i}$, and $M_{\chi,i}$. So they are
available for us to use here.  Our goal here is to redo the work in
Section~\ref{sec:same-orbit} without having a $\mathbf{0''}$ listing
of the sets $\Dhat_\chi$, $\Rhat_{\chi,i}$, and $\Mhat_{\chi,i}$.  Our
goal is to find a suitable listing of these sets and the isomorphisms
$\Phi_{\chi,\Lambda(\chi)}$.  And then start working from
Section~\ref{sec:build-phiI} onward to construct the desired
automorphism using the replacement parts we have constructed.  We work
with an oracle for $\Lambda$ and $0^{(8)}$.

$\Lambda$ is an isomorphism between $T$ and $\mathcal{T}(\Ahat)$.  By
Lemma~\ref{oracle}, using $\mathbf{0}^{(8)}$ as an oracle, we can find
a representative of each atomic $T$-equivalence class of
$n$-successors of $\Ahat$. Furthermore, we can assume that when
choosing a representative we always choose a maximal representative
of terms of $T$-equivalence.  Hence we can consider $\Lambda$ as a map
that takes $D_\chi$ to a representative of the equivalent class which
codes $\chi$.  Let $\Dhat_{\Lambda(\chi)} = \Lambda(D_\chi)$.

We recall that each $R_{\chi,i}$ is broken into a number of pieces.
First there is a subset $M_{\chi,i}$ which is either maximal in
$R_{\chi,i}$ or almost equal to $R_{\chi,i}$.  $M_{\chi,i}$ is split
into several parts; $H_{\chi,i}$ and if $\xi = \chi\widehat{~}l \in T$
and $l^{-1}(\xi,0) \leq l^{-1}(\chi,i)$ or $\xi = \chi$ then $D_\xi
\cap M_{\chi,i} = D_{\xi} \cap R_{\chi,i}$ is a infinite split of
$M_{\chi,i}$; $M_{\chi,i}$ is computable iff all of these pieces are
computable.  Effectively in each $\chi$ and $i$ we can give a finite
set $F_{\chi,i}$ such that
$$ R_{\chi,i} = (R_{\chi,i} - M_{\chi,i}) \sqcup 
H_{\chi,i} \sqcup \bigsqcup_{\xi \in F_{\chi,i}} (D_{\xi} \cap
R_{\chi,i})$$ and either, for all $\xi \in F_{\chi,i}$, $M_{\chi,i}$
is maximal in $R_{\chi,i}$ and $D_\xi \cap R_{\chi,i}$ is a nontrivial
split of $M_{\chi,i}$ or, for all $\xi \in F_{\chi,i}$, $M_{\chi,i} =
R_{\chi,i}$ and $D_\xi \cap R_{\chi,i}$ is computable.  Now we must
find $\Rhat_{\Lambda(\chi),i}$ such that it has the same properties.

We need the following two lemmas.  The first follows from the
definition of an extendible subalgebra.  The second lemma follows from
the construction of $A_T$ and the fact that, for almost all $i$,
$D_{\xi}$ lives in $R_{\xi^-,i}$ iff $D_{\xi^-}$ lives in
$R_{\xi^-,i}$.  The second part of the second lemma follows in
particular from the homogeneity requirements.

\begin{lemma}
  The collection of the sets
  \begin{equation}
    \begin{split}
      \{(R_{\xi^-,i} \cap D_\xi) : i \geq j\},
      \{(\overline{R}_{\xi^-,i} \cap D_\xi): i \geq j\}, \\
      \{(R_{\xi,i} \cap D_\xi):i \geq 0 \}, \text{ and } \{
      (\overline{R}_{\xi,i} \cap D_\xi) : i\geq 0\}
    \end{split}
  \end{equation}
  form an extendible subalgebra, $\mathbb{B}_{\xi,j}$, of the splits
  of $D_\xi$.
\end{lemma}

\begin{lemma}\label{same}
  If $|\xi| = |\zeta|$ then there is a $j_{\xi,\zeta}$ such that
  $\mathbb{B}_{\xi,j_{\xi,\zeta}}$ is extendibly
  $\Delta^0_3$-isomorphic to $\mathbb{B}_{\zeta,j_{\xi,\zeta}}$ via
  the identity map. (The identity map sends $R_{\xi,i} \cap D_\xi$ to
  $R_{\zeta,i} \cap D_\zeta$, etc.)  Furthermore, for all $i$, $D_\xi$
  lives in $R_{\chi, i}$ iff $D_\zeta$ lives in $R_{\chi,i}$ and, for
  all $i \geq j_{\xi,\zeta}$, $D_\xi$ lives in $R_{\chi^-, i}$ iff
  $D_\zeta$ lives in $R_{\chi^-,i}$. 
\end{lemma}

Now we must use another theorem from \citet{Cholak.Harrington:nd}.


\begin{theorem}[Theorem~6.3 of \citet{Cholak.Harrington:nd}]\label{old2}
  Assume $D$ and $\Dhat$ are automorphic via $\Psi$. Then $D$ and
  $\Dhat$ are automorphic via $\Theta$ where $\Theta \upharpoonright
  \E(D)$ is $\Delta^0_3$. 
\end{theorem}

\begin{lemma}\label{onlyDelta03}
  For some $j_\xi$, there is an extendible subalgebra,
  $\widehat{\mathbb{B}}_{\Lambda(\xi),j_\xi}$, of the splits of
  $D_{\Lambda(\xi)}$ which is extendibly $\Delta^0_3$ isomorphic via
  $\Theta_\xi$ to $\mathbb{B}_{\xi,j_\xi}$.  Furthermore, for all $i
  \geq j_\xi$, $D_\xi \cap R_{\xi^-,i}$ is the split of a maximal set
  iff $\Theta_\xi(D_\xi \cap R_{\xi^-,i})$ is the split of a maximal
  set, and $D_\xi \cap R_{\xi^-,i}$ is computable iff
  $\Theta_\xi(D_\xi \cap R_{\xi^-,i})$ is computable. And, for all
  $i$, $D_\xi \cap R_{\xi,i}$ is the split of a maximal set iff
  $\Theta_\xi(D_\xi \cap R_{\xi,i})$ is the split of a maximal set,
  and $D_\xi \cap R_{\xi,i}$ is computable iff $\Theta_\xi(D_\xi \cap
  R_{\xi,i})$ is computable.
  Moreover, we can find $j_\xi$,
  $\widehat{\mathbb{B}}_{\Lambda(\xi),j_\xi}$, and $\Theta_\xi$ with
  an oracle for $\mathbf{0}^{(8)}$.
\end{lemma}

\begin{proof}
  Recall $A$ and $\Ahat$ are automorphic via $\Psi$ and the image of a
  $D_\xi$ must also code a node of length $|\xi|$. By
  Lemma~\ref{codingstocodings}, $\Dhat_{\Lambda(\xi)}$ is the pre-image
  under $\Psi$ of some $D_{\Psi^{-1}(\Lambda(\xi))} =^* D_\eta -
  \bigsqcup_{j < j'} R_{\eta^-,j}$, where $|\eta| = |\xi|$. Now apply
  Theorem~\ref{old2} to get $\Theta_\xi$.  Find the least $j_\xi$ such
  that, for all $i \geq j_\xi$, $D_{\Lambda(\xi)}$ lives in
  $R_{\Lambda(\xi)^-,i}$ iff $D_{\Lambda(\xi)^-}$ lives in
  $R_{\Lambda(\xi)^-,i}$ and similarly for
  $D_{\Psi^{-1}(\Lambda(\xi))}$ and $D_{\Psi^{-1}(\Lambda(\xi))^-}$,
  and $D_\xi$ and $D_{\xi^-}$.  The image of
  $\mathbb{B}_{\Psi^{-1}(\Lambda(\xi)),j_\xi}$ under $\Theta_\xi$ is
  an extendible subalgebra $\widehat{\mathbb{B}}_{\Lambda(\xi),j_\xi}$
  and, furthermore, these subalgebras are extendibly
  $\Delta^0_3$-isomorphic.  By Lemma~\ref{same},
  $\mathbb{B}_{\xi,j_\xi}$ is extendibly $\Delta^0_3$-isomorphic to
  $\mathbb{B}_{\Psi^{-1}(\Lambda(\xi)),j_\xi}$. Since $\Theta_\xi$ is
  an automorphism the needed homogeneous properties are preserved.

  Now that we know these items exist we know that we can successfully
  search for them.  Look for an $j_\xi$ and $\Theta_\xi$ such that
  $\Theta_\xi(\mathbb{B}_{\xi,j_\xi}) =
  \widehat{\mathbb{B}}_{\Lambda(\xi),j_\xi}$ is extendibly
  $\Delta^0_3$-isomorphic to $\mathbb{B}_{\xi,j_\xi}$ via
  $\Theta_\xi$; these items also satisfiy the second sentence of the
  above lemma and the additional property that, for all $\Rhat$, if
  $\Rhat$ is an infinite subset of $D_{\Lambda(\xi)}$ then there are
  finitely many $\Rtilde_i$ such that $\Rhat \subseteq^* \bigcup
  \Theta_\xi(\Rtilde_i)$.  Since, by Requirement~\ref{sec:requireI},
  this last property is true of $D_\xi$, and $\Theta_\xi$ is generated
  by an automorphism, it also must be true of $D_{\Lambda(\xi)}$.
  This extra property ensures that $\Theta_\xi$ is onto.  By carefully
  counting quantifiers we see that $\mathbf{0}^{(8)}$ is more than
  enough to find these items.
\end{proof}

Let $\tilde{F}_{\chi,i}$ be such that $\xi \in \tilde{F}_{\chi,i}$ iff
$\xi \in {F}_{\chi,i}$ and $i \geq j_\xi$.  For all $\chi$ and $i$,
let $$\breve{\Hhat}_{\Lambda(\chi),i} = \bigsqcup_{\xi \in
  \tilde{F}_{\chi,i}} \Theta_\xi(D_{\xi} \cap R_{\chi,i}).$$ Either
$\breve{\Hhat}_{\Lambda(\chi),i}$ is computable or the split of a
maximal set. This follows from the projection through the above lemmas
of the homogeneity requirements.  In the latter case,
$\breve{\Hhat}_{\Lambda(\chi),i}$ lives inside $\widehat{\omega}$.

We repeatly apply the dual of Lemma~\ref{sec:livesI} to all those
$\breve{\Hhat}_{\Lambda(\chi),i}$ who live inside $\widehat{\omega}$
to get $\Rtilde_{\Lambda(\chi),i}$ which are all pairwise disjoint.
This determines the $\tilde{\Mhat}_{\Lambda(\chi),i}$ which witness
that $\breve{\Hhat}_{\Lambda(\chi),i}$ lives in
$\Rtilde_{\Lambda(\chi),i}$.  Let $\breve{\Rhat}_{\Lambda(\chi),i}$ be
a computable infinite subset of $\tilde{\Mhat}_{\Lambda(\chi),i} -
\breve{\Hhat}_{\Lambda(\chi),i}$ (we call this \emph{set
  subtraction}).  Let $\Rhat_{\Lambda(\chi),i} =
\Rtilde_{\Lambda(\chi),i} - \breve{\Rhat}_{\Lambda(\chi),i}$.
$\breve{\Hhat}_{\Lambda(\chi),i}$ lives inside
$\Rhat_{\Lambda(\chi),i}$.  In this case, again, by the dual of 
Lemma~\ref{sec:livesI}, we have determined $\Mhat_{\Lambda(\chi),i}$
and hence we have determined $\widehat{H}_{\Lambda(\chi),i}$.

So it remains to find $\Rhat_{\Lambda(\chi),i}$ and
$\Mhat_{\Lambda(\chi),i}$, where $\breve{\Hhat}_{\Lambda(\chi),i}$ is
computable.  For such $i$ once we find $\Rhat_{\Lambda(\chi),i}$ we
will let $\Rhat_{\Lambda(\chi),i} = \Mhat_{\Lambda(\chi),i}$.

By Requirement~\ref{sec:requireI} and our construction, for all $e$,
there are finite sets $F_D$ and $F_R$ such that either
\begin{equation*}
  W_e \subseteq^* \biggl ( \bigsqcup_{\chi \in F_D}D_\chi \cup
  \bigsqcup_{(\chi,i) \in F_R} R_{\chi,i} \biggr ),
\end{equation*}
or there is an $R_{W_e}$ such that
\begin{equation*}
  R_{W_e} \subseteq \bigg ( \bigsqcup_{\chi \in F_D}D_\chi \cup
  \bigsqcup_{(\chi,i) \in F_R} R_{\chi,i} \bigg) \text{ and } W_e \cup
  R_{W_e} = \omega.
\end{equation*}
By Lemma~\ref{codingstocodings}, as a collection the
$\Dhat_{\Lambda(\chi)}$s are the isomorphic images of the collection
of the $D_\chi$ and similarly with the collection of all
$R_{\chi,i}$s. 
 Hence we should be able to define
$\Rhat_{\Lambda(\chi),i}$, where $\breve{\Hhat}_{\Lambda(\chi),i}$ is
computable such that, for all $e$, there are finite sets $\Fhat_D$ and
$\Fhat_R$ with either
\begin{equation} \label{eq:dualeq1} \What_e \subseteq^* \biggl (
  \bigsqcup_{\chi \in \Fhat_D} \Dhat_{\Lambda(\chi)} \cup
  \bigsqcup_{(\chi,i) \in \Fhat_R} \Rhat_{\Lambda(\chi),i} \biggr ),
\end{equation}
or there is an $R_{\What_e}$ such that
\begin{equation}\label{eq:dualeq2}
  R_{\What_e} \subseteq \bigg ( \bigsqcup_{\chi \in
    \Fhat_D}\Dhat_{\Lambda(\chi)} \cup \bigsqcup_{(\chi,i) \in
    \Fhat_R} \Rhat_{\Lambda(\chi),i} \bigg) \text{ and } \What_e \cup
  R_{\What_e} = \widehat{\omega}.
\end{equation}

Fix some nice listing of the $(\chi,i)$ such that
$\Rhat_{\Lambda(\chi),i}$ has yet to be defined (as above).  Assume
that $(\chi,i)$ is the $e$th member in our list and the first $e-1$ of
$\Rhat_{\Lambda(\chi),i}$ have been defined such that, for all $e' <
e$, one of the two equations above hold.  For all $e$, either there
are finitely many $(\xi,j)$ where $\Rhat_{\Lambda(\xi),j}$ is defined such
that $\Rhat_{\Lambda(\xi),j} \cap \What_e \neq^* \emptyset$ or, for
almost all $(\xi,j)$, where $\Rhat_{\Lambda(\xi),j}$ is defined,
$\Rhat_{\Lambda(\xi),i} \subseteq^* \What_e$ (this is true for any
possible pre-image of $\What_e$).

In the first case find a computable $\Rhat$, a finite $\Fhat_R$, and a
finite $\Fhat_D$ such that if $(\xi,j) \in \Fhat_R$ then
$\Rhat_{\Lambda(\xi),j}$ is defined; if $\Rhat_{\Lambda(\xi),j}$ is
defined then $\Rhat \cap \Rhat_{\Lambda(\xi),j} = \emptyset$;
$\breve{\Hhat}_{\Lambda(\chi),i} \subseteq \Rhat$; $(\Rhat -
\breve{\Hhat}_{\Lambda(\chi),i}) \cap \bigsqcup_\xi
\Dhat_{\Lambda(\xi)} = \emptyset$ (these last three clauses are
possible because of the above set subtraction); and \[\What_e
\subseteq^* \biggl ( \Rhat \cup \bigsqcup_{\xi \in \Fhat_D}
\Dhat_{\Lambda(\xi)} \cup \bigsqcup_{(\xi,i) \in \Fhat_R}
\Rhat_{\Lambda(\xi),i} \biggr ).\]

In the second case find a computable $\Rhat$, a finite $\Fhat_R$, and
a finite $\Fhat_D$ such that all of the above but the last clause
above hold and
\[\overline{\What_e} \subseteq^* \biggl ( \Rhat \cup \bigsqcup_{\xi
  \in \Fhat_D} \Dhat_{\Lambda(\xi)} \cup \bigsqcup_{(\xi,i) \in
  \Fhat_R} \Rhat_{\Lambda(\xi),i} \biggr ).\] Either way let
$R_{\Lambda(\chi),i} = \Rhat$.  Since the sets we have defined so far
cannot be all the images of the $R_{\xi,l}$, there must be enough of
$\widehat{\omega}$ for us to continue the induction.



Now we have to find a replacement for the isomorphisms given to us by
Lemma~\ref{isomorphismnice}; we cannot.  But as we work through
Section~\ref{sec:build-phiI} we see that we want to apply Theorem~5.10
of \citet{Cholak.Harrington:nd} to $D_{\xi} - \bigsqcup_{j<j_\xi}
R_{\xi^-,j}$ and $D_{\Lambda(\xi)} - \bigsqcup_{j<j_\xi}
\Theta_\xi(R_{\xi^-,j})$, we need these isomorphisms to meet the
hypothesis, and, furthermore, this is the only place these isomorphisms
are used. However, the first step of the proof of Theorem~5.10 of
\citet{Cholak.Harrington:nd} is to use the given isomorphisms (given by
Lemma~\ref{isomorphismnice}) to create an extendible isomorphism
between extendible subalgebra generated by $R_{\chi,i} \cap D_\xi$ and
the one generated by $\Rhat_{\chi,i} \cap \Dhat_{\Lambda(\xi)}$ and,
furthermore, this is the only place these given isomorphisms are used
in the proof. These subalgebras are $\mathbb{B}_{\xi,j_\xi}$ and
$\widehat{\mathbb{B}}_{\Lambda(\xi),j_\xi}$ which are isomorphic via
$\Theta_\xi$.  Hence we can assume that we can apply Theorem~5.10 of
\citet{Cholak.Harrington:nd}.

At this point we have all the needed sets and isomorphisms with the
desired homogeneity between these sets (in terms of
Requirement~\ref{sec:homo}).  Now we have enough to apply part 3 of
our game plan to construct the desired automorphism.  That is, start
working from Section~\ref{sec:build-phiI} onward to construct the
desired automorphism.

\subsection{A Technical Invariant for the orbit of $A_T$}\label{tinv}
The goal of this section is to prove a theorem like
Theorem~~\ref{sec:invariant} but without the hypothesis that $A$ and
$\Ahat$ are in the same orbit.  Reflecting back through the past
section we see that the fact that $A$ and $\Ahat$ are in the same
orbit was used twice: in the proof of Lemma~\ref{onlyDelta03} and in
showing that Equations~\ref{eq:dualeq1} and \ref{eq:dualeq2} hold.
Hence we assume these two items would allow us to weaken the hypothesis
as desired. Since the notation from the above section is independent
of the fact that $A$ and $\Ahat$ are in the same orbit we borrow it
wholesale for the following.

\begin{theorem}
  Assume
  \begin{enumerate}
  \item $T \cong \mathcal{T}(\Ahat)$ via $\Lambda$,
  \item the conclusion of Lemma~\ref{onlyDelta03} (the whole statement
    of the lemma is the conclusion), and
  \item Equations~\ref{eq:dualeq1} and \ref{eq:dualeq2} hold.
  \end{enumerate} 
  Then $A_T \approx \Ahat$ via $\Phi_\Lambda$ where
  $\Phi_\Lambda\leq_T \Lambda \join \bf{0}^{(8)}$.
\end{theorem}

\begin{corollary}\label{techsameorbit}
  $A_T \approx \Ahat$ iff \begin{enumerate}
  \item $T \cong \mathcal{T}(\Ahat)$ via $\Lambda$,
  \item the conclusion of Lemma~\ref{onlyDelta03} (the whole statement
    of the lemma is the conclusion), and
  \item Equations~\ref{eq:dualeq1} and \ref{eq:dualeq2} hold.
  \end{enumerate}
\end{corollary}

\section{Our Orbits and Hemimaximal Degrees}\label{hemimaximal}

A set is \emph{hemimaximal} iff it is the nontrivial split of a
maximal set.  A degree is \emph{hemimaximal} iff it contains a
hemimaximal set.

Let $T$ be given.  Construction $A_T$ as above.  For all $i$, either
$A_T$ lives in $R_i$ or $A_T \cap R_i$ is computable. If $A_T$ lives
in $R_i$ then $A_T \cap R_i$ is a split of maximal set $M \sqcup
\overline{ R}_i$ and hence $A_T = (A_T \cap R_i)$ is a hemimaximal
set.  $A_T = \bigsqcup_{i \in \omega} (A_T \cap R_i)$ where $A_T \cap
R_i$ is either hemimaximal or computable.  So the degree of $A_T$ is
the infinite join of hemimaximal degrees.  It is not known if the
(infinite) join of hemimaximal degrees is hemimaximal.  Moreover, this
is not an effective infinite join. But if we control the degrees of
$A_T \cap R_i$ we can control the degree of $A_T$.

\begin{theorem}\label{AThemi}
  Let $H$ be hemimaximal.  We can construct $A_T$ such that $A_T
  \equiv_T H$.  Call this $A_T$, $A^H_T$, to be careful.
\end{theorem}

\begin{proof}
  Consider those $\alpha$ and $k$ such that $l(|\alpha|-k) =
  (\lambda,n)$, for some $n$.  Only at such $\alpha$ do we construct
  pieces of $D^k_\lambda = A_{T_k}$.  Uniformly we can find partial
  computable mapping, $p^k_\alpha$, from $\omega$ to $R^k_\alpha$ such
  that if $R^k_\alpha$ is an infinite computable set then $p^k_\alpha$
  is one-to-one, onto, and computable.  Since $H$ is hemimaximal there
  is a maximal set $M$ and a split $\breve{H}$ witnessing that $H$ is
  hemimaximal.  Then $ p^k_\alpha(M) \sqcup \overline{R}^k_\alpha$ is
  maximal and $p^k_\alpha(H)$ is nontrivial split of $ p^k_\alpha(M)
  \sqcup \overline{R}^k_\alpha$ with the same degree as $H$.

  The idea is that at $\alpha$ we would like to let $M^k_\alpha =
  p^k_\alpha(M)$ but because of the dumping this does not work.
  Dumping allows us to control whether $R^k_\alpha =^* M^k_\alpha$ or
  not.  Let $\tilde{M}^k_\alpha = p^k_\alpha(M)$. If
  $$\overline{p^k_\alpha(M_s)} \cap {R}^k_\alpha = \{
  m^{\alpha,k}_0,m^{\alpha,k}_1,m^{\alpha,k}_2, \ldots\}$$ then place
  the marker $\Gamma^{\alpha,k}_e$ on $m^{\alpha,k}_e$ at stage $s$.
  Now when dumping the element marked by marker $\Gamma^{\alpha,k}_e$
  we will just dump that \emph{single} element (this not the case in
  the standard dumping arguments).  Now assume that the dumping is
  done effectively (this is the case in the construction of $A_T$).
  Let $M^k_{\alpha,s+1} = \tilde{M}^k_{\alpha,s+1} \cup
  M^k_{\alpha,s}$ plus those $m^\alpha_e$ which are dumped via
  $\Gamma^{\alpha,k}_e$ at stage $s+1$.  $M^k_\alpha$ is c.e.\ and
  $\tilde{M}^k_\alpha \subseteq M^k_\alpha$.  Since
  $\tilde{M^k_\alpha} \sqcup \overline{R}^k_\alpha$ is maximal, either
  $M^k_\alpha =^* \tilde{M}^k_\alpha$ or $M^k_\alpha =^*
  R^k_\alpha$. In the first case $p^k_\alpha(H)$ and
  $p^k_\alpha(\breve{H}) \sqcup \overline{R}^k_\alpha$ are nontrivial
  splits of $M^k_\alpha$.  The second case occurs iff there is least
  $\Gamma^{\alpha,k}_e$ which is dumped into $M^k_\alpha$ infinitely
  often.  The above construction of $M^k_\alpha$ is uniformly in
  $\alpha$.

  In Section~\ref{sec:maximal-sets-their}, when we construct
  $M^k_\alpha$ and its splits, rather than using the maximal set
  construction and the Friedberg splitting construction, we use the
  above construction of $M^k_\alpha$; we will put the split
  $p^k_\alpha(H)$ into $D^k_\lambda = A_T$ and use the Friedberg
  splitting construction to split $p^k_\alpha(\breve{H})$ into enough
  pieces as determined by the construction.
\end{proof}

There is no reason to believe that if $\Ahat$ is in the same orbit as
$A^H_T$ that $\Ahat \equiv_T H$. Nor is there a reason to believe $\Ahat$
must have hemimaximal degree.  Notice that for each $H$ we have a
separate construction. Hence the homogeneity requirement need not hold
between these different constructions.  Therefore, we cannot prove
that the sets $A^H_T$ are in the same orbit. It might be that for $H
\neq \tilde{H}$, that $A^H_T$ and $A^{\tilde{H}}_T$ are in different
orbits.  We conjecture, using Corollary~\ref{techsameorbit}, it is
possible to construct two different versions of $A_T$ which are not in
the same orbit.  But we can do the following.

\begin{theorem}
  There is an $A_T$ whose orbits contain a representative of every
  hemimaximal degree.
\end{theorem}

\begin{proof}
  The idea is for all hemimaximal $H$ to do the above construction
  simultaneously.  This way the homogeneous requirement will be met
  between the different $A^H_T$s.

  Notice the above construction is uniformly in the triple $e =
  \langle m,h,\breve{h}\rangle$ where $W_m = M, W_h = H$, and
  $W_{\breve{h}} = \breve{H}$.
  
  We want to reorder the trees from Theorem~\ref{folkloreI}. Let
  $\tilde{T}_{\langle e,i \rangle} = T_i$.  Now do the construction in
  Section~\ref{mainproof} with two expectations: use the trees
  $\tilde{T}_{\langle e,i \rangle}$ and, for those $\alpha$ and $k$
  such that $l(|\alpha|-k) = (\lambda,n)$, for some $n$, we use the
  construction of $M^k_\alpha$ outlined in the proof of
  Theorem~\ref{AThemi}.
   
  For all $i$ and $e$ coding a hemimaximal set we construct a set
  $A_{\tilde{T}_{\langle e, i \rangle}}$.  If $e'$ codes another
  hemimaximal set then $A_{\tilde{T}_{\langle e, i \rangle}}$ and
  $A_{\tilde{T}_{\langle e', i \rangle}}$ are in the same orbit.

  If $e'$ does not code sets such that $W_m = W_h \sqcup
  W_{\breve{h}}$ then construction of $A_{\tilde{T}_{\langle e', i
      \rangle}}$ is impaired but this does not impact the simultaneous
  construction of the other $A_{\tilde{T}_{\langle e, i \rangle}}$.
\end{proof}

\section{On the Isomorphism Problem for Boolean Algebras and Trees}

\nocite{Sacks:90} \nocite{Rogers:67}

\subsection{$\Sigma^1_1$-completeness}\label{folkloreproofI}

We think it is well known that the isomorphism problem for Boolean
Algebras and Trees are $\Sigma^1_1$-complete, at least in the form
stated in Theorems~\ref{folkloreBA} and \ref{folkloreI}. We have
searched for a reference to a proof for these theorems without
success.  It seems very likely that these theorems were known to
Kleene. There are a number of places where something close to what we
want appears; for example, see \citet{walkerwhite:00},
\citet{MR2033315}, and the example at the end of Section~5 of
\citet{MR2058190}. Surely there are other examples.  All of these work
by coding the Harrison ordering, as will the construction below.  To
be complete we include a proof in this section.  The material we
present below is similar to results in the three papers mentioned
above. We are thankful to Noam Greenberg for providing the included
proof.

\begin{remark}[Notation]
  For cardinals $\kappa,\lambda$, etc.\ (we use $2$ and $\w$), a
  \emph{tree on $\kappa\times \lambda$} is a downward-closed subset of
  \[ \bigcup_{n<\w} \kappa^n \times \lambda^n ,\] so that the set of
  paths of the tree is a closed subset of $\kappa^\w\times
  \lambda^\w$. We may use more or fewer coordinates. For a tree $R$,
  $[R]$ is the set of paths through $R$. For a subset $A$ of a product
  space $\kappa^\omega \times \lambda^\w$ (for example), $pA$ is the
  projection of $A$ onto the first coordinate.
\end{remark}

\begin{lemma}
  There is an effective operation $I$ such that, given a computable
  infinite-branching tree $T$, $I(T)$ is a computable linear ordering
  such that
  \begin{enumerate}
  \item if $T$ is well-founded then $I(T)$ is a well-ordering;
  \item if $T$ is not well-founded then $I(T) \cong \wock(1+\Q)$.
  \end{enumerate}
\end{lemma}

\begin{proof}
  Suppose that a computable tree $T_0\subseteq \w^{<\w}$ is given.
  Unpair to get a tree $T_1$ on $2\times \w$ such that $[T_0] =
  \{X\oplus f\,:\, (X,f)\in [T_1]\}.$

  Now let $T_2 = T_1 \times 2^{<\w}$, the latter inserted as a second
  coordinate (so $T_2 = \{ (\sigma,\tau,\rho)\,:\, (\sigma,\rho)\in
  T_1 \andd \tau\in 2^{<\w} \andd |\tau| = |\sigma| = |\rho|\}$.) Let
  $T_3$ be the tree on $2\times \w$ which is obtained by pairing the
  first two coordinates of $T_2$.

  The class $\HYP$ of hyperarithmetic reals is $\Pi^1_1$, and so
  $p[T_3] - \HYP$ is $\Sigma^1_1$; let $T_4$ be a computable tree such
  that $p[T_4] = p[T_3] - \HYP$.

  Let $L_5$ be the Kleene-Brouwer linear ordering obtained from $T_4$;
  finally, let $I(T) = L_5\w = L_5 + L_5 + \cdots$.

  The point is this: $p[T_2] = p[T_1]\times 2^\w$. Thus if $T$ is not
  well-founded, then $p[T_1]$ is nonempty and so $p[T_2]$ is
  uncountable and so $p[T_4]$, and hence $[T_4]$, is nonempty. If $T$
  is well-founded then $p[T_4]$ is empty; that is, $T_4$ is
  well-founded. Also, $p[T_4]$ contains no hyperarithmetic sets, and
  so $T_4$ has no hyperarithmetic paths.

  It follows that if $T$ is well-founded then $L_5$, and so $I(T)$, is
  a well-ordering. If $T$ is not well-founded then $L_5$ is a
  computable linear ordering which is not a well-ordering but has no
  hyperarithmetic infinite descending chains, that is, a Harrison linear
  ordering. This has order-type $\wock(1+\Q) + \gamma$ for some
  computable ordinal $\gamma$. For any computable $\gamma$ we have
  $\gamma + \wock = \wock$ (as $\wock$ is closed under addition) and
  so $I(T)$ has ordertype $\wock(1+\Q+1+\Q+1+\Q+\cdots) \cong
  \wock(1+\Q)$.
\end{proof}

\begin{corollary}[Proposition~5.4.1 of \citet{walkerwhite:00}]
  For any $\Sigma^1_1$ set $A$, there is a computable sequence
  $\seq{L_n}$ of (computable) linear orderings such that, for all $n$,
  \begin{enumerate}
  \item if $n\in A$ then $L_n \cong \wock(1+\Q)$;
  \item if $n\notin A$ then $L_n$ is a well-ordering.
  \end{enumerate}
\end{corollary}

\begin{proof}
  Let $A$ be a $\Sigma^1_1$ set. There is a computable sequence
  $\seq{T_n}$ of trees on $\w$ such that, for all $n$, $n\notin A$ iff
  $T_n$ is well-founded.  Now apply $I$ to each $T_n$.
\end{proof}

\begin{corollary}[Theorem~\ref{folkloreI}]\label{cor32}
  There is a computable tree $T$ on $\w$ such that the collection of
  computable trees $S$ which are isomorphic to $T$ is
  $\Sigma^1_1$-complete.
\end{corollary}

\begin{proof}
  Use the operation that converts a linear ordering $L$ to the tree
  $T_L$ of finite descending sequences in $L$. The point is that if
  $L$ is an ordinal then $T_L$ is well-founded and so cannot be
  isomorphic to $T_{\wock(1+\Q)}$.
\end{proof}

\begin{corollary}[Theorem~\ref{folkloreBA}]
  There is a computable Boolean algebra $B$ such that the collection
  of Boolean algebras $C$ that are isomorphic to $B$ is
  $\Sigma^1_1$-complete.
\end{corollary}

\begin{proof}
  Similar; use the interval algebra $B_L$. If $L$ is an ordinal then
  $B_L$ is superatomic.
\end{proof}

\subsection{$\Pi^0_n$-completeness}\label{folkloreproofII}

Again we believe it is known that there are trees $T_{\Pi_n}$ such
that the isomorphism problem for $T_{\Pi_n}$ is $\Pi^0_n$-complete, at
least in the form stated in Theorem~\ref{sec:folkloredark}.  The
closest we could find was work in \citet{walkerwhite:00}, which does
not quite work.  To be complete we include a proof in this
section. The details are similar in style but different from what
is found in \cite{walkerwhite:00}. The trees in \cite{walkerwhite:00}
 do not provide precise bounds; they are hard for the appropriate
class but not known to be complete (see Remark~\ref{walker}). We
wonder if Theorem~\ref{sec:folkloredark} is true for all computable
ordinals, the case $\alpha=\omega$ being a good test case.  The
following construction is joint work with Noam Greenberg.  The
following lemma is well known, but we include a proof for
completeness; it is a partial version of uniformalization.

\begin{lemma}
  Let $A(n,x)$ be a $\Pi^0_1$ relation. Then there is a $\Pi^0_1$
  partial function $f$ such that $\dom A = \dom f$.
  \label{lem: one witness}
\end{lemma}

\begin{proof}
  We give an effective construction of a computable predicate $R$ such
  that $f(n)=x \iff \forall y\, R(n,x,y)$. If $n\ge s$ or $x\ge s$
  then $R(n,x,s)$ always holds; so to make $R$ computable, at stage
  $s$ of the construction we define $R(n,x,s)$ for all $x,n<s$. In
  fact, for all $n<s$, at stage $s$ we define $R(n,x,s)$ to hold for
  at most one $x<s$. This will imply that $f$ is indeed a function.

  Let $S$ be a computable predicate such that $A(n,x) \iff \forall y\,
  S(n,x,y)$.

  For every $n$ and $x$ we have a moving marker $c(n,x)$. We start
  with $c(n,x)=x$. At stage $s$, for every $n<s$, find the least $x<s$
  such that for all $y<s$ we have $S(n,x,y)$ (if one exists). For
  $x'\ne x$, initialize $c(n,x')$ by redefining it to be large. Now
  define $R$ by letting $R(n,c(n,x),s)$ hold but $R(n,z,s)$ not hold
  for all $z<s$ different from $c(n,x)$.

  Let $n<\w$. Suppose that $n\in \dom f$. For all $s>\max\{n,f(n)\}$,
  $R(n,f(n),s)$ holds, which means that at stage $s$, $f(n) = c(n,x)$
  for some $x$.  Different markers get different values and so there
  is just one such $x$, independent of $s$. By the instructions, for
  all $s>\max\{n,f(n)\}$, for all $y<s$, $S(n,x,y)$ holds; this shows
  that $n\in \dom A$.

  Suppose that $n\in \dom A$. Let $x$ be the least such that for all
  $y$, $S(n,x,y)$ holds. There is some stage after which $c(n,x)$ does
  not get initialized (wait for some stage $s$ that bounds, for all
  $z<x$, some $y$ such that $S(n,z,y)$ does not hold). Let $s$ be the
  last stage at which $c(n,x)$ gets initialized. At stage $s$, a
  final, large value $a=c(n,x)$ is chosen. For all $t>a$, $R(n,a,t)$
  holds because $t>s$. Thus $a$ witnesses that $n\in \dom f$.
\end{proof}

By relativizing the above to $\cero^{(n-2)}$, we see that for every
$n\ge 2$, for every $\Sigma^0_n$ set $A$, there is a $\Pi^0_{n-1}$
function $f$ such that $A = \dom f$.

A \emph{tree} is a downward closed subset of $\w^{<\w}$. The
collection $\Tree$ of all computable trees (i.e., indices for total,
computable characteristic functions of trees) is $\Pi^0_2$. For any
tree $T$, let $\Isom_T$ be the collection of $S\in \Tree$ which are
isomorphic to $T$.

\begin{lemma}
  Let $T_{\Pi_2}$ be the infinite tree of height 1.
  $\Isom_{T_{\Pi_2}}$ is $\Pi^0_2$-complete.
\end{lemma}

\begin{proof}
  A tree is isomorphic to $T_{\Pi_2}$ iff it has height 1 and it is
  infinite.  Certainly this is a $\Pi^0_2$ property.

  Let $A$ be a $\Pi^0_2$ set; say that $A(n) \iff \forall x \exists y
  R(n,x,y)$ where $R$ is computable. For $n$ and $s$, let $l(n,s)$ be
  the greatest $l$ such that for all $x\le l$ there is some $y<s$ such
  that $R(n,x,y)$ holds. Say that $s$ is expansionary for $n$ if
  $l(n,s) > l(n,s-1)$.

  For each $n$ define a tree $T_{2,A}(n)$: this is a tree of height 1,
  and a string $\seq{s}$ is on the tree iff $s$ is expansionary for
  $n$. Then $n\mapsto T_{2,A}(n)$ reduces $A$ to $\Isom_{T_{\Pi_2}}$.
\end{proof}

For the next level we use trees of height 2. We use two trees: the
tree $T_{\Pi_3}$ is the tree of height 2 such that for each $n$ there
are infinitely many level 1 nodes which have exactly $n$ children, and
no level 1 node has infinitely many children. The tree $T_{\Sigma_3}$
is like $T_{\Pi_3}$, except that we add one level 1 node which has
infinitely many children.

\begin{lemma}
  $\Isom_{T_{\Pi_3}}$ is $\Pi^0_3$ and $\Isom_{T_{\Sigma_3}}$ is
  $\Pi^0_3 \land \Sigma^0_3$.
\end{lemma}

\begin{proof}
  If $T$ is a computable tree, then the predicate ``$\seq{x}$ has
  exactly $n$ children in $T$'' is $\Sigma^0_2$, uniformly in a
  computable index for $T$.  So is the predicate ``$\seq{x}$ has
  finitely many children in $T$''. The predicate ``there are
  infinitely many level 1 nodes on $T$ which have $n$ children'' is
  $\Pi^0_3$.

  Also, to say that the height of a tree $T$ is at most 2 is $\Pi^0_1$
  (once we know that $T\in \Tree$).

  A tree $T$ is isomorphic to $T_{\Pi_3}$ if it has height at most 2
  and for every $n$, there are infinitely many level 1 nodes on $T$
  which have $n$ children, and every level 1 node on $T$ has finitely
  many successors.

  The predicate ``$\seq{x}$ has infinitely many children in $T$'' is
  $\Pi^0_2$; and so the predicate ``at most one level 1 node on $T$
  has infinitely many children'' is $\Pi^0_3$.

  A tree $T$ is isomorphic to $T_{\Sigma_3}$ if it has height at most
  2 and for every $n$, there are infinitely many level 1 nodes on $T$
  which have $n$ children, at most one level 1 node on $T$ has
  infinitely many children, and some level 1 node has infinitely many
  children. The last condition is $\Sigma^0_3$ and all previous ones
  are $\Pi^0_3$.
\end{proof}

\begin{lemma} $(\Sigma^0_3,\Pi^0_3) \le_1 (\Isom_{T_{\Sigma_3}},
  \Isom_{T_{\Pi_3}})$.
\end{lemma}

\begin{proof}
  Let $A$ be a $\Sigma^0_3$ set. By Lemma~\ref{lem: one witness},
  there is some $\Pi^0_2$-definable function $f$ such that $A = \dom
  f$.

  For any $n$, we define a tree $T_{3,A}(n)$ of height 2. First, it
  contains a copy of $T_{\Pi_3}$. Then, for every $x$, there is a
  level 1 node $\seq{m_x}$ such that $T_{3,A}(n)[m_x] = T_{2,f}(n,x)$
  (that is, for all $y$, $\seq{m_x,y} \in T_{3,A}(n)$ iff $\seq{y}\in
  T_{2,f}(n,x)$.)

  Then $n\mapsto T_{3,A}(n)$ reduces $(A,\lnot A)$ to
  $(\Isom_{T_{\Sigma_3}}, \Isom_{T_{\Pi_3}})$ because for all but
  perhaps one $x$ we have $T_{2,f}(n,x)$ finite.
\end{proof}

\begin{remark}[Walker's $T_{\Sigma_3}$]\label{walker}
  Walker defined his $T_{\Sigma_3}$ such that it has infinitely many
  $T_{\Pi_2}$ children. Walker's $\Isom_{T_{\Sigma_3}}$ is be
  $\Pi^0_4$.  The above lemma still holds (via a slightly different
  reduction) but we only get hardness not completeness.  It is not
  known if Walker's $T_{\Sigma_3}$ is $\Pi_4$-complete.  To avoid
  using infinitely many $T_{\Pi_2}$ children we have to be more
  careful.  Here we get around this problem by using
  Lemma~\ref{lem: one witness}.
\end{remark}

We can now lift it up.

\begin{lemma} For all $n\ge 3$ there are trees $T_{\Sigma_n}$ and
  $T_{\Pi_n}$ such that
  \begin{enumerate}
  \item $\Isom_{T_{\Pi_n}}$ is $\Pi^0_n$;
  \item $\Isom_{T_{\Sigma_n}}$ is $\Pi^0_n \land \Sigma^0_n$;
  \item $(\Sigma^0_n, \Pi^0_n) \le_1 (\Isom_{T_{\Sigma_n}},
    \Isom_{T_{\Pi_n}})$.
  \end{enumerate}
\end{lemma}

Thus $\Isom_{T_{\Pi_n}}$ is $\Pi^0_n$-complete.

\begin{proof} By induction; we know this for $n=3$.

  The tree $T_{\Pi_{n+1}}$ is a tree of height $n$ which has
  infinitely many level 1 nodes, the tree above each of which is
  $T_{\Sigma_n}$. The tree $T_{\Sigma_{n+1}}$ is the tree
  $T_{\Pi_{n+1}}$, together with one other level 1 node above which we
  have $T_{\Pi_n}$.

  A tree $T$ is isomorphic to $T_{\Pi_{n+1}}$ iff it has infinitely
  many level 1 nodes (this is $\Pi^0_2$!), and for every level 1 node
  $\seq{x}$, the tree $T[x]$ above $\seq{x}$ is isomorphic to
  $T_{\Sigma_n}$.

  A tree $T$ is isomorphic to $T_{\Sigma_{n+1}}$ iff it has infinitely
  many level 1 nodes; for every level 1 node $\seq{x}$, the tree
  $T[x]$ is isomorphic to either $T_{\Sigma_n}$ or to $T_{\Pi_n}$;
  there is at most one $\seq{x}$ such that $T[x]$ is isomorphic to
  $T_{\Pi_n}$; and there is some $\seq{x}\in T$ such that $T[x] \cong
  T_{\Pi_n}$.

  Note again that if we had infinitely many $T_{\Pi_n}$s (which is
  what White's trees had) then we'd have had to pay another
  quantifier.

  The reduction is similar to that of the case $n=3$: given a
  $\Sigma^0_{n+1}$ set $A$, we get a $\Pi^0_n$ function $f$ such that
  $A = \dom f$; we construct $T_{n+1,A}(m)$ to be a tree such that for
  all $x$, $\seq{x}\in T_{n+1,A}(m)$ and the tree $T_{n+1,A}(m)[x] =
  T_{n,f}(m,x)$.
\end{proof}

For the case $\alpha\geq \omega$, the situation is murkier. Using the
trees from \citet{walkerwhite:00}, for example, gives a reduction of,
say, $\Sigma^0_{\w+1}$ to a tree $T$ such that $\Isom_{T}$ is
computable from something like $\cero^{(\w+3)}$. With more work it
seems that this can be reduced to $\cero^{(\w+2)}$, but it seems
difficult to reduce this to $\cero^{(\w)}$.  We remark that ``things
catch up with themselves'' at limit levels which is why we get $+2$
for $\alpha\geq \omega$.


\begin{thebibliography}{20}
\providecommand{\natexlab}[1]{#1}
\providecommand{\url}[1]{\texttt{#1}}
\expandafter\ifx\csname urlstyle\endcsname\relax
  \providecommand{\doi}[1]{doi: #1}\else
  \providecommand{\doi}{doi: \begingroup \urlstyle{rm}\Url}\fi

\bibitem[Ash and Knight(2000)]{MR1767842}
C.~J. Ash and J.~Knight.
\newblock \emph{Computable structures and the hyperarithmetical hierarchy},
  volume 144 of \emph{Studies in Logic and the Foundations of Mathematics}.
\newblock North-Holland Publishing Co., Amsterdam, 2000.
\newblock ISBN 0-444-50072-3.

\bibitem[Cholak(2006)]{2006:c}
Peter Cholak.
\newblock The {C}omputably {E}numerable {S}ets: the {P}ast, the {P}resent and
  the {F}uture.
\newblock Theory and Applications of Models of Computation, 2006, Beijing
  China, 2006.

\bibitem[Cholak and Harrington(2005)]{Cholak.Harrington:nd}
Peter Cholak and Leo~A. Harrington.
\newblock Extension theorems, orbits, and automorphisms of the computably
  enumerable sets.
\newblock To appear in Trans. Amer. Math. Soc. Final version as of 8/31/2005.
  math.LO/0408279, 2005.

\bibitem[Cholak and Harrington(2002)]{mr2003h:03063}
Peter Cholak and Leo~A. Harrington.
\newblock On the definability of the double jump in the computably enumerable
  sets.
\newblock \emph{J. Math. Log.}, 2\penalty0 (2):\penalty0 261--296, 2002.
\newblock ISSN 0219-0613.

\bibitem[Cholak and Harrington(2003)]{mr2004f:03077}
Peter Cholak and Leo~A. Harrington.
\newblock Isomorphisms of splits of computably enumerable sets.
\newblock \emph{J. Symbolic Logic}, 68\penalty0 (3):\penalty0 1044--1064, 2003.
\newblock ISSN 0022-4812.

\bibitem[Cholak et~al.(2001)Cholak, Downey, and Herrmann]{mr2001k:03086}
Peter Cholak, Rod Downey, and Eberhard Herrmann.
\newblock Some orbits for {$\mathcal{E}$}.
\newblock \emph{Ann. Pure Appl. Logic}, 107\penalty0 (1-3):\penalty0 193--226,
  2001.
\newblock ISSN 0168-0072.

\bibitem[Downey and Stob(1992)]{Downey.Stob:92}
R.~G. Downey and M.~Stob.
\newblock Automorphisms of the lattice of recursively enumerable sets: Orbits.
\newblock \emph{Adv. in Math.}, 92:\penalty0 237--265, 1992.

\bibitem[Goncharov et~al.(2004)Goncharov, Harizanov, Knight, and
  Shore]{MR2058190}
Sergey~S. Goncharov, Valentina~S. Harizanov, Julia~F. Knight, and Richard~A.
  Shore.
\newblock {$\Pi^1_1$} relations and paths through {$\mathcal{O}$}.
\newblock \emph{J. Symbolic Logic}, 69\penalty0 (2):\penalty0 585--611, 2004.
\newblock ISSN 0022-4812.

\bibitem[Harrington and Soare(1998)]{Harrington.Soare:98}
Leo Harrington and Robert~I. Soare.
\newblock Codable sets and orbits of computably enumerable sets.
\newblock \emph{J. Symbolic Logic}, 63\penalty0 (1):\penalty0 1--28, 1998.
\newblock ISSN 0022-4812.

\bibitem[Harrington and Soare(1991)]{Harrington.Soare:91}
Leo~A. Harrington and Robert~I. Soare.
\newblock {P}ost's program and incomplete recursively enumerable sets.
\newblock \emph{Proc. Nat. Acad. Sci. U.S.A.}, 88:\penalty0 10242--10246, 1991.

\bibitem[Harrington and Soare(1996)]{Harrington.Soare:96}
Leo~A. Harrington and Robert~I. Soare.
\newblock The ${\Delta}\sp 0\sb 3$-automorphism method and noninvariant classes
  of degrees.
\newblock \emph{J. Amer. Math. Soc.}, 9\penalty0 (3):\penalty0 617--666, 1996.
\newblock ISSN 0894-0347.

\bibitem[Hirschfeldt and White(2002)]{MR2033315}
Denis~R. Hirschfeldt and Walker~M. White.
\newblock Realizing levels of the hyperarithmetic hierarchy as degree spectra
  of relations on computable structures.
\newblock \emph{Notre Dame J. Formal Logic}, 43\penalty0 (1):\penalty0 51--64
  (2003), 2002.
\newblock ISSN 0029-4527.

\bibitem[Lachlan(1968)]{Lachlan:68*3}
Alistair~H. Lachlan.
\newblock On the lattice of recursively enumerable sets.
\newblock \emph{Trans. Amer. Math. Soc.}, 130:\penalty0 1--37, 1968.

\bibitem[Maass(1984)]{Maass:84}
W.~Maass.
\newblock On the orbit of hyperhypersimple sets.
\newblock \emph{J. Symbolic Logic}, 49:\penalty0 51--62, 1984.

\bibitem[Rogers(1967)]{Rogers:67}
H.~Rogers, Jr.
\newblock \emph{Theory of Recursive Functions and Effective Computability}.
\newblock McGraw-Hill, New York, 1967.

\bibitem[Sacks(1990)]{Sacks:90}
Gerald Sacks.
\newblock \emph{Higher Recursion Theory}.
\newblock Perspectives in Mathematical Logic. Springer--Verlag, Heidelberg,
  1990.

\bibitem[Slaman and Woodin(1989)]{Slaman.Woodin:conjecture}
Theodore~A. Slaman and W.~Hugh Woodin.
\newblock Slaman-{W}oodin conjecture.
\newblock Personal communication, 1989.

\bibitem[Soare(1974)]{Soare:74}
Robert~I. Soare.
\newblock Automorphisms of the lattice of recursively enumerable sets {I}:
  maximal sets.
\newblock \emph{Ann. of Math. (2)}, 100:\penalty0 80--120, 1974.

\bibitem[Soare(1987)]{Soare:87}
Robert~I. Soare.
\newblock \emph{Recursively Enumerable Sets and Degrees}.
\newblock Perspectives in Mathematical Logic, Omega Series. Springer--Verlag,
  Heidelberg, 1987.

\bibitem[White(2000)]{walkerwhite:00}
Walker~M. White.
\newblock \emph{Characterizations for Computable Structures}.
\newblock PhD thesis, Cornell University, Ithaca, NY, USA, 2000.

\end{thebibliography}
\end{document}